\pdfoutput=1
\documentclass[10pt,reqno]{amsart}

\numberwithin{equation}{section}

\usepackage[margin=1in]{geometry}
\usepackage{cancel}
\usepackage[tt=false]{libertine}
\usepackage{mathtools}
\usepackage{amssymb}
\usepackage{mathrsfs}
\usepackage[varbb]{newpxmath}

\let\savedbigtimes\bigtimes
\let\bigtimes\relax
\usepackage{mathabx} 
\let\bigtimes\savedbigtimes

\usepackage{textcomp}
\usepackage{enumerate}
\usepackage{bm}

\usepackage[usenames,dvipsnames]{xcolor}
\usepackage[colorlinks=true,
	linkcolor=teal!60!black,
	citecolor=PineGreen,
	urlcolor=RedViolet]{hyperref}
\hypersetup{bookmarksopen=true}
\usepackage[footnotesize,
	justification=centering]{caption}
\usepackage[font=footnotesize]{subcaption}
\usepackage{graphicx}

\newtheorem{thm}{Theorem}[section]
\newtheorem{lem}[thm]{Lemma}
\newtheorem{ppn}[thm]{Proposition}
\newtheorem{cor}[thm]{Corollary}

\theoremstyle{definition}
\newtheorem{dfn}[thm]{Definition}

\newtheorem*{rmk*}{Remark}
\newtheorem{ass}[thm]{Assumption}

\definecolor{Green}{HTML}{239B56}
\definecolor{Blue}{HTML}{2471A3}
\definecolor{Orange}{HTML}{D35400}
\definecolor{DarkBlue}{HTML}{1B4F72}
\definecolor{DarkGray}{HTML}{1C2833}

\newcommand{\beq}{\begin{equation}}
\newcommand{\eeq}{\end{equation}}

\newcommand{\<}{\langle}
\renewcommand{\>}{\rangle}

\newcommand{\f}{\frac}
\newcommand{\set}[1]{\{#1\}}
\newcommand{\Ind}[1]{\mathbf{1}\{#1\}}
\newcommand{\I}{\mathbf{1}}

\DeclareMathOperator{\Var}{Var}

\DeclareMathOperator{\dist}{dist}

\newcommand{\E}{\mathbb{E}}
\renewcommand{\P}{\mathbb{P}}
\newcommand{\R}{\mathbb{R}}

\title[Sharp threshold sequence and universality for Ising perceptron models]{Sharp threshold sequence and universality for Ising perceptron models}
\author[S.\ Nakajima]{Shuta Nakajima$^\star$}
\author[N.\ Sun]{Nike Sun$^\circ$}
\date{\today}
\thanks{$^\star$Graduate School of Science and Technology, Meiji University, Tokyo. $^\circ$Department of Mathematics, Massachusetts Institute of Technology.}

\newcommand{\RR}{\bar{R}}
\newcommand{\rr}{\bar{r}}
\newcommand{\ZETA}{\zeta} 
\newcommand{\XI}{\xi} 
\newcommand{\bXI}{\bm{\xi}} 
\newcommand{\NU}{\nu} 

\newcommand{\filt}{\mathscr{F}}
\newcommand{\GG}{\mathscr{G}}
\newcommand{\HH}{\mathscr{H}}

\newcommand{\wmax}{w_{\max}}
\newcommand{\bDe}{\bm{\Delta}}
\newcommand{\mm}{\bm{m}}
\newcommand{\bS}{\bm{S}}

\begin{document}

\maketitle

\begin{abstract} We study a family of Ising perceptron models with $\{0,1\}$-valued activation functions. This includes the classical half-space models, as well as some of the symmetric models considered in recent works. For each of these models we show that the free energy is self-averaging, there is a sharp threshold sequence, and the free energy is universal with respect to the disorder. A prior work of C.\ Xu (2019) used very different methods to show a sharp threshold sequence in the half-space Ising perceptron with Bernoulli disorder. Recent works of Perkins--Xu (2021) and Abbe--Li--Sly (2021) determined the sharp threshold and limiting free energy in a symmetric perceptron model. The results of this paper apply in more general settings, and are based on new ``add one constraint'' estimates extending Talagrand's estimates for the half-space model (1999, 2011).
\end{abstract}

\setcounter{tocdepth}{1}
\tableofcontents

\section{Introduction}
\label{s:intro}

The \textbf{perceptron} is a classical model in high-dimensional probability theory \cite{cover1965geometrical,MR146858}, which can be interpreted as a toy model of a simple neural network \cite{gardner1988space}. In the simplest version, the \textbf{Ising (half-space) perceptron} refers to the random set $\bS$ defined by intersecting the discrete cube $\set{-1,+1}^N$ with $M=N\alpha$ i.i.d.\ random half-spaces (formal definitions below). In the late 1980s it was shown that heuristic analytical tools of statistical physics (i.e., the ``replica method'' or ``cavity method'') can be applied to derive precise predictions on the limiting behavior of the perceptron model \cite{gardner1988space,gd1988optimal,krauthmezard1989,mezard1989space}. In particular, it was predicted that the \textbf{free energy} $N^{-1}\log|\bS|$ of the model concentrates around an explicit constant (depending on the parameters of the model), and that there is a \textbf{sharp threshold} $\alpha_c$ such that $\P(|\bS|>0)$ transitions from $1-o_N(1)$ to $o_N(1)$ in an $o_N(1)$ window around $\alpha_c$.

In contrast with what has been conjectured via physics heuristics, mathematical understanding of Ising perceptron models remains quite limited. For the half-space Ising perceptron, it follows by a trivial first moment calculation that $\P(|\bS|>0) = o_N(1)$ for $\alpha=M/N$ large enough, but this bound is typically not tight. On the other hand, it was shown  by \cite{MR1629627} (see also \cite{MR1716771}) that  $\P(|\bS|>0) = 1-o_N(1)$ for $\alpha=M/N$ small enough. Moreover, for small enough $\alpha$, the free energy concentrates around the value conjectured by physicists \cite{MR1782273,MR3024566,bolthausen2021gardner}. A result of \cite{ding2018capacity} shows that, under an additional condition, the random set $\bS$ is nonempty with \emph{positive} probability
for any $\alpha$ smaller than the critical value predicted by physicists \cite{krauthmezard1989}. These results leave open the question of whether the model has a sharp threshold, and whether the free energy concentrates for general $\alpha$. 

For the half-space Ising perceptron with \textbf{Bernoulli disorder} --- meaning that the half-spaces point in directions chosen uniformly at random from $\set{-1,+1}^N$ --- C.\ Xu applied Hatami's pseudo-junta theorem \cite{MR2925389} to show that the model has a \textbf{sharp threshold sequence} \cite{MR4317708}. In a different direction, numerous recent works have studied \textbf{symmetric} perceptron models, which are significantly more tractable than the half-space models and can be analyzed for much finer properties \cite{MR3983947,perkins2021frozen,abbe2021proof,abbe2021binary,gamarnik2022algorithms}. In particular, for a \textbf{symmetric interval} variant of the Ising perceptron (intersecting $\set{-1,+1}^N$ with i.i.d.\ random symmetric \emph{slabs} rather than half-spaces), recent work rigorously pinpoints the limiting free energy and sharp threshold \cite{perkins2021frozen,abbe2021proof}.

In this paper we study a family of Ising perceptron models with $\set{0,1}$-valued activation functions (defined more formally below). This includes the classical half-space models as well as some of the symmetric models considered in more recent works. Further, we allow the disorder random variables to come from any \textbf{subgaussian} distribution, in contrast with previous works which have assumed gaussian or Bernoulli disorder. For each of these models we show:
\begin{itemize}
\item There is a \textbf{sharp threshold sequence}, meaning that
 $\P(|\bS|>0)$ transitions from $1-o_N(1)$ to $o_N(1)$ in an $o_N(1)$ window of $\alpha$ values, although the location of the transition may depend on $N$.
\item The model is \textbf{self-averaging}, meaning that the free energy $N^{-1}\log|\bS|$ concentrates around a deterministic value which may depend on $N$.
\item The free energy and sharp threshold sequence are \textbf{universal} with respect to the disorder.  
\end{itemize}
Our sharp threshold sequence result extends the main theorem of \cite{MR4317708}, although by very different methods --- our approach does not use Hatami's theorem or other tools of boolean analysis. 

For the half-space Ising perceptron model with gaussian or Bernoulli disorder, the sharp threshold sequence and self-averaging property appear to be relatively straightforward consequences of ``add one constraint'' estimates proved in \cite{MR1680236,MR3024566} (and discussed below). While Talagrand does indicate the implications of his estimates for self-averaging, to the best of our knowledge he does not explicitly address the sharp threshold problem; and one of the aims of this article is to clearly spell out this simple connection. The main technical contributions of this paper are new ``add one constraint'' estimates for general (not necessarily symmetric) interval perceptron models estimates with general (subgaussian) disorder. These results cannot be deduced by an easy generalization of Talagrand's arguments, and we describe some of the main ideas in \S\ref{ss:intro.perceptron} below. We then apply our ``add one constraint'' estimates to deduce the sharp threshold sequence and free energy concentration results. Lastly we combine these with some central limit theorem considerations to deduce the universality results. A weaker version of our universality results (comparing gaussian with Bernoulli disorder, and with additional smoothness assumptions) appears in \cite[\S9.9]{MR3024566}.

The remainder of this introductory section is organized as follows:
\begin{itemize}
\item In \S\ref{ss:intro.perceptron} we formally define the various perceptron models that we consider. We state our first set of results (Theorems~\ref{t:hspace.general}--\ref{t:interval.general}) which control
the effect of adding a single constraint in these perceptron models. We also describe some of the basic proof ideas.

\item In \S\ref{ss:conc.results} we state our second set of results (Theorems~\ref{t:conc}--\ref{t:univ})  which address self-averaging, sharp thresholds, and universality in the perceptron models. These may be viewed as consequences of the ``add one constraint'' estimates described in \S\ref{ss:intro.perceptron}.

\item In \S\ref{ss:intro.lit} we discuss the context given by the surrounding literature on perceptron models.
\end{itemize}
At the end of this section we give the outline for the rest of the paper.

\subsection{``Add one constraint'' estimates for perceptron models}
\label{ss:intro.perceptron}

In the above we used $\bS$ to denote the Ising perceptron solution set, a random subset of $\set{-1,+1}^N$. In this paper we only consider the size of the random set, and we hereafter always denote $Z\equiv |\bS|$. We will prove results for the classical Ising perceptron model 
	\beq\label{e:Z.hspace}
	Z \equiv Z_{M,N}(\kappa,\infty,\XI)
	\equiv \sum_{\sigma\in\set{-1,+1}^N}
		\prod_{k\le M} \I\bigg\{
		\f{(\XI^k,\sigma)}{N^{1/2}} \ge \kappa
		\bigg\}\,,
	\eeq
which we will refer to as the \textbf{half-space (Ising) perceptron}. We also consider more general models of the form
	\beq\label{e:Z.U}
	Z_{M,N}(U;\XI)
	\equiv 
	\sum_{\sigma\in\set{-1,+1}^N}
	\prod_{k\le M}
	U\bigg( \f{(\XI^k,\sigma)}{N^{1/2}} \bigg)
	\eeq
where $U:\R\to[0,1]$ (measurable) is the \textbf{activation function}. We refer to this as the $U$-perceptron; most of this paper concerns the case that $U$ is $\set{0,1}$-valued. We refer to the case $U(x)=\Ind{x\in[a,b]}$ as the \textbf{interval (Ising) perceptron}.
The \textbf{symmetric interval perceptron} mentioned above corresponds to the case $a=-b$. We assume that the $\XI, \XI^k$ ($k\ge1$) are i.i.d.\ random vectors satisfying the following:

\begin{ass}\label{a:subgaus}
The $\XI_i$ are i.i.d.\ random variables with mean zero, unit variance, such that
	\[
	\E\exp(\lambda\XI_i)
	\le \exp\bigg(\f{\lambda^2\NU}{2}\bigg)\]
for all $\lambda\in\R$ (that is, the $\XI_i$ are subgaussian with variance proxy $\NU$). This assumption includes the two most widely studied cases of the perceptron model: \textbf{gaussian disorder} (the $\XI^k$ are standard gaussian vectors) and \textbf{Bernoulli disorder} (the $\XI^k$ are sampled uniformly at random from $\set{-1,+1}^N$). Note that in general we must have $\NU\ge1$, since the cumulant-generating function $\mathfrak{K}(\lambda)\equiv\log\E \exp(\lambda\XI_i)$ satisfies 
$\mathfrak{K}''(0)=\Var\XI_i=1$.
\end{ass}

We first restate a result of Talagrand:

\begin{thm}[half-space perceptron, gaussian disorder \cite{MR1680236,MR3024566}]\label{t:hspace.gaus}
Let $S$ be any subset of $\set{-1,+1}^N$ with $|S|\ge \exp(N\delta)$, where $\delta$ is a small positive constant. If $g$ is a standard gaussian random vector in $\R^N$, then
	\[\P\bigg(
	\f1{|S|}
		\bigg|\bigg\{\sigma\in S: 
	\f{(g,\sigma)}{N^{1/2}}
		\ge \kappa 
		\bigg\}\bigg|
	< \f1{\exp(w)}
	\bigg)
	\le \exp\bigg( -\f{w}{C_\delta}\bigg)\,.
	\]
provided $C_\delta \le w \le N/C_\delta$, where $C_\delta$ is a large finite constant that depends only on $\kappa$ and $\delta$.
\end{thm}

Theorem~\ref{t:hspace.gaus} was in fact already proved by Talagrand in the more difficult setting of Bernoulli disorder (see \cite[Propn.~2.3]{MR1680236}). The result for the gaussian case follows by a simplification of the argument for the Bernoulli case, which we review in \S\ref{ss:talagrand}. Most of the argument for the gaussian case also appears in \cite[Ch.~9]{MR3024566}. (The latter reference \cite{MR3024566} does not appear to treat the case of small $\delta$, but this was already done previously in \cite{MR1680236}.)

As we discuss further in \S\ref{ss:sup.conc}--\ref{ss:sudakov} below,
the proof of Theorem~\ref{t:hspace.gaus} relies on two main ingredients:
\begin{enumerate}[(i)]
\item a lower bound on the expected supremum of a stochastic process; and
\item a concentration estimate for the supremum of a stochastic process.
\end{enumerate}
In the gaussian setting these ingredients are supplied by the well-known \textbf{Sudakov minoration lower bound} and \textbf{Borell--TIS (Tsirelson--Ibragimov--Sudakov) concentration inequality}, restated in Lemmas~\ref{l:sudakov} and \ref{l:borell.tis} below. For more general distributions, however, both of these ingredients can be quite nontrivial. 

By contrast, our next result gives a weaker bound than Theorem~\ref{t:hspace.gaus}, but has the advantage that it uses only the simplest gaussian versions of the results mentioned above (Sudakov minoration and Borell--TIS). We transfer these bounds to general distributions by a rather weak form of the multivariate CLT. As a result, the bound of Theorem~\ref{t:hspace.general} is almost certainly suboptimal, but it has the advantage that it is relatively straightforward to derive, and applies for a general class of subgaussian distributions. To the best of our knowledge, using existing methods, the argument of Talagrand for Theorem~\ref{t:hspace.gaus} may only be extended to a limited subclass of subgaussian distributions (see \S\ref{ss:sup.conc}--\ref{ss:sudakov} for the details).

\begin{thm}[half-space perceptron, general disorder]\label{t:hspace.general}
Let $S$ be any subset of $\set{-1,+1}^N$ with $|S|\ge \exp(N\delta)$, where $\delta$ is a small positive constant. Suppose $\XI$ is a random vector in $\R^N$ satisfying Assumption~\ref{a:subgaus}.  Then
	\[\P\bigg(
	\f1{|S|}
		\bigg|\bigg\{\sigma\in S: 
	\f{(\XI,\sigma)}{N^{1/2}}
		\ge \kappa 
		\bigg\}\bigg|
	< \f1{\exp(w)} \bigg)
	\le \exp\bigg( -\f{w}{C_\delta}\bigg)\,.\]
provided $C_\delta \le w \le N^{1/2}/C_\delta$, where $C_\delta$ is a large finite constant that depends only on $\kappa$, $\NU$, and $\delta$.\end{thm}

The next two results concern the $U$-perceptron \eqref{e:Z.U}, under the assumption that we have $U:\R\to\set{0,1}$ with $U(x)\ge\Ind{x\in[a,b]}$ for some $-\infty<a<b<\infty$. 

\begin{thm}[$U$-perceptron, gaussian disorder]\label{t:interval.gaus}
Let $S$ be any subset of $\set{-1,+1}^N$ with $|S|\ge\exp(N\delta)$ where $\delta$ is a small positive constant. If $g$ is a standard gaussian random vector in $\R^N$, then
	\[\P\bigg( 
	\f1{|S|}
	\bigg|\bigg\{\sigma\in S: 
		\f{(g,\sigma)}{N^{1/2}}
		\in[a,b]\bigg\}\bigg|
	< \f1{\exp(w)}
	\bigg) \le \exp\bigg( -\f{w}{C_\delta}\bigg)
	\]
provided $C_\delta\le w \le N^{1/2}/C_\delta$, where $C_\delta$ is a large finite constant that depends only on $a$, $b$, and $\delta$.
\end{thm}

\begin{thm}[$U$-perceptron, general disorder]
\label{t:interval.general}
Let $S$ be any subset of $\set{-1,+1}^N$ with $|S|\ge\exp(N\delta)$ where $\delta$ is a small positive constant. Suppose $\XI$ is a random vector in $\R^N$ satisfying Assumption~\ref{a:subgaus}. Then
	\[\P\bigg( 
	\f1{|S|}
	\bigg|\bigg\{\sigma\in S: 
		\f{(\XI,\sigma)}{N^{1/2}}
		\in[a,b]\bigg\}\bigg|
	< \f1{\exp(w)}
	\bigg) \le \exp\bigg(- \f{w}{C_\delta}\bigg)
	\]
provided $C_\delta\le w\le N^{1/3}/C_\delta$, where $C_\delta$ is a large finite constant that depends only on $a$, $b$, $\NU$, and $\delta$.
\end{thm}

We now describe some of the basic ideas that appear in the proofs of Theorems~\ref{t:hspace.general}--\ref{t:interval.general}. We repeatedly leverage a basic observation that comes from Talagrand's proof of Theorem~\ref{t:hspace.gaus}: if $S$ is a large subset of $\set{-1,+1}^N$, then it must contain many elements that are well separated in Hamming distance. If $g$ is a standard gaussian vector in $\R^N$, then the separation guarantee can be used to show that
	\[
	\max\bigg\{ \f{(g,\sigma)}{N^{1/2}} : \sigma\in S\bigg\}
	\]
is likely to be quite large. The details of this argument are reviewed in \S\ref{ss:talagrand} below. We point out that it relies on theorems for gaussian processes, namely, concentration of the supremum and a Sudakov minoration lower bound on the expectation of the supremum. As we discuss in \S\ref{ss:sup.conc}--\ref{ss:sudakov}, these do not easily extend to more general subgaussian distributions.

In this paper we devise a more flexible version of Talagrand's approach, as follows. First we decompose $[N]\equiv\set{1,\ldots,N}$ into $L$ blocks of size $K$ each, where $KL=N$. Denote the blocks $I_1,\ldots,I_L$. For any subset $I\subseteq [N]$ let $\sigma_I\equiv (\sigma_i)_{i\in I}$; we will say that two configurations $\sigma,\tau\in\set{-1,+1}^N$ are well-separated on $I$ if $\sigma_I$ and $\tau_I$ are not too close (see Definition~\ref{d:sep} below). We show that if $S\subseteq\set{-1,+1}^N$ is large, then it must contain many elements that are well-separated on a positive fraction of the blocks --- say, on the final $L\gamma$ blocks. We then consider the process
	\beq\label{e:block.process}
	M_k(\sigma)
	= \f1{N^{1/2}}
		\sum_{j=1}^k (\XI_{I_J},
			\sigma_{I_j})
	\eeq
for $\sigma\in S$, starting from $M_0(\sigma)=0$ and with the goal of having say $M_L(\sigma)\in[a,b]$ (for the $U$-perceptron \eqref{e:Z.U}, assuming $U(x)\ge\Ind{x\in[a,b]}$). In the case that $\XI$ is \textbf{gaussian}, we can apply Talagrand's estimates for the half-space model \eqref{e:Z.hspace} to control the increments of $M_k(\sigma)$ over the final $L\gamma$ steps, using the separation guarantees on the final $L\gamma$ blocks. With this approach we show that, with very good probability, $S$ must contain many elements $\sigma$ with $M_L(\sigma)\in[a,b]$.

In the case that $\XI$ follows a more general \textbf{subgaussian} distribution, we follow a similar strategy of first obtaining estimates for the half-space model \eqref{e:Z.hspace}, then using a block decomposition to obtain estimates for the $U$-perceptron \eqref{e:Z.U}. However, as we already noted, it does not appear that Talagrand's results for the half-space perceptron extend easily to the full class of subgaussian distributions. Instead, to prove the desired estimates in the half-space model, we again use the block decomposition approach, combined with a (rather weak) form of the multivariate central limit theorem which allows us to borrow the estimates from the gaussian setting. In particular, we restrict ourselves to quantities of the form
	\[
	\max\bigg\{ \f{(\XI,\sigma)}{N^{1/2}} : \sigma\in X\bigg\}
	\] 
for subsets $X\subseteq\set{-1,+1}^N$ of \emph{bounded} size, allowing us to apply central limit theorems only in bounded dimensions. This results in weaker bounds which can potentially be improved by appealing to stronger forms of the central limit theorem. 

\subsection{Self-averaging, sharp threshold sequence, and universality} 
\label{ss:conc.results}

We next describe the main results which we obtain as consequences of the ``add one constraint'' estimates presented in \S\ref{ss:intro.perceptron}. To emphasize the dependence on the ``add one constraint'' estimates, we state some of our results for a more abstract model of the form
	\beq\label{e:Z.Theta}
	Z \equiv Z_{M,N} \equiv \sum_{\sigma\in\set{-1,+1}^N}
	\prod_{k\le M} \Theta_k(\sigma)\,,
	\eeq
where the $\Theta_k$ are a sequence of $[0,1]$-valued random functions,  adapted to a filtration $(\filt_k)_{k\ge0}$, satisfying a fairly weak ``add one constraint'' estimate:

\begin{ass}\label{a:addone.weak}
For the model \eqref{e:Z.Theta}, let $Z_{M+1}$ be the partition function that results from introducing one more factor $\Theta_{M+1}$. Suppose that for all $\delta>0$ small enough that on the event $Z_M \ge \exp(N\delta)$ we have
	\[
	\P\bigg(
	\f{Z_{M+1}}{Z_M} \le \f1{\exp(w)} \,\bigg|\,
	\filt_M\bigg) \le f_\delta(w)\]
for all $C_\delta \le w \le \wmax$, such that the function $f_\delta$ satisfies the bound
	\beq\label{e:assumed.tail.bound}
	\int_{C_\delta}^{\wmax} 
	2w f_\delta(w)\,dw
	\le C_{\delta,2}\,.
	\eeq
In the above, $C_\delta$ and $C_{\delta,2}$ are finite constants that depend only on the model and on $\delta$, while $\wmax$ can depend on the model as well as on $\delta$ and $N$.
\end{ass}

The half-space and $U$-perceptron models \eqref{e:Z.hspace} and \eqref{e:Z.U} are clearly a special case of the model \eqref{e:Z.Theta}. Therefore, Theorems~\ref{t:hspace.gaus}--\ref{t:interval.general} imply that the perceptron models \eqref{e:Z.hspace} and \eqref{e:Z.U} 
satisfy a  stronger condition than Assumption~\ref{a:addone.weak}; see Assumption~\ref{a:addone} at the start of Section~\ref{s:sharp}. 
However, we will show that the weaker Assumption~\ref{a:addone.weak} suffices for some similar results. Denote $\log_{N\delta}Z\equiv \max\set{\log Z,N\delta}$.\footnote{This differs slightly from the notation of \cite[\S8.3]{MR3024566} which has
$\log_A x = \max\set{\log x,-A}$. We choose a notation which is more convenient for our choice of normalization.}

\begin{thm}[concentration of free energy]
\label{t:conc}
For the model \eqref{e:Z.Theta}, under Assumption~\ref{a:addone.weak}, we have
\[\lim_{N\to\infty}
	\bigg| \log_{N\delta}Z-\E\log_{N\delta}Z\bigg| =0\]
for any positive constant $\delta$,
where the limit holds in probability as $N\to\infty$.
\end{thm}

Our next result concerns the behavior of the probability $\P(Z_M>0)$ for the model \eqref{e:Z.Theta}. For this question, we may as well assume that the $\Theta_a$ are $\set{0,1}$-valued, since otherwise we can replace $\Theta_a(x)$ with $\Ind{\Theta_a(x)>0}$. We introduce one more assumption below, which together with Assumption~\ref{a:addone.weak} ensures that  $\P(Z_M>0)$ has a transition in the proportional regime $M\asymp N$. This assumption is typically easily verifiable in practice: 

\begin{ass}\label{a:firstmmtubd}
There is a positive constant $c$ such that
$\E( \Theta_{k+1}\,|\,\filt_k)\le \exp(-c)$ almost surely for all $k\ge0$.
\end{ass}

\begin{thm}[sharp threshold sequence]\label{t:sharpthreshold}
For the model \eqref{e:Z.Theta}, if the $\Theta_a$ are $\set{0,1}$-valued and satisfy Assumptions~\ref{a:addone.weak} and \ref{a:firstmmtubd}, then there is a sharp threshold sequence: that is to say, there is a sequence $\alpha_N\asymp1$ such that
$\P(Z_{N\alpha,N}>0)$ transitions from $1-o_N(1)$ to $o_N(1)$ in an $o_N(1)$ window around $\alpha_N$.\end{thm}

In particular, by Theorems~\ref{t:hspace.gaus}--\ref{t:interval.general}, 
the results of Theorems~\ref{t:conc} and \ref{t:sharpthreshold} apply to the perceptron models \eqref{e:Z.hspace} and \eqref{e:Z.U}, assuming in the latter model that $U$ is $\set{0,1}$-valued with $U(x)\ge\Ind{x\in[a,b]}$, and further assuming the $\XI^k$ satisfy Assumption~\ref{a:subgaus}. For the sharp threshold question, we can replace $U(x)$ by $\bar{U}(x)=\Ind{U(x)>0}$, so we can allow the original $U$ to take all values in $[0,\infty]$. Thus Theorem~\ref{t:sharpthreshold} implies the following: 

\begin{cor}\label{c:sharpthreshold}
For a measurable function $U:\R\to[0,1]$, let $\bar{U}(x)=\Ind{U(x)>0}$. If $\bar{U}(x)\ge\Ind{x\in[a,b]}$ for some $-\infty<a<b<\infty$,
and $\bar{U}$ is not almost everywhere one, then the results of Theorem~\ref{t:sharpthreshold} hold for the model \eqref{e:Z.U}, provided that the $\XI^k$ satisfy Assumption~\ref{a:subgaus}.

\begin{proof}
Since we are only interested in whether the partition function $Z$ of \eqref{e:Z.U} is nonnegative, we may assume without loss that $U=\bar{U}$. By the condition $U(x)\ge\Ind{x\in[a,b]}$, combined with Assumption~\ref{a:subgaus} and Theorem~\ref{t:interval.general}, the model \eqref{e:Z.U} satisfies Assumption~\ref{a:addone}. The condition that $U$ is not almost everywhere one guarantees that Assumption~\ref{a:firstmmtubd} is also satisfied.
The claim then directly follows from Theorem~\ref{t:sharpthreshold}.
\end{proof}
\end{cor}

Finally, we can combine the above results with some central limit theorem estimates to obtain the following:

\begin{thm}[universality]\label{t:univ}
In the perceptron model \eqref{e:Z.U}, suppose $U$ is $\set{0,1}$-valued and piecewise continuous, and neither identically one nor identically zero. Suppose the $\XI^k$ are i.i.d.\ random vectors satisfing Assumption~\ref{a:subgaus}. Then the threshold sequence and free energy are universal with respect to the disorder.
\end{thm}

\subsection{Related work}\label{ss:intro.lit}
As mentioned earlier, the perceptron model was analyzed in the (nonrigorous) statistical physics literature in the 1980s \cite{gardner1988space,gd1988optimal,krauthmezard1989,mezard1989space}. We refer to \cite{cover1965geometrical,gardner1988space,montanari2019generalization,montanari2021tractability} for discussions of statistical motivations for the model. 

There is an important variant of the perceptron model which is well understood at a rigorous level:  the \textbf{spherical (half-space) perceptron}, defined by intersecting the sphere $N^{1/2}\mathbb{S}^{N-1}$ with i.i.d.\ random \textbf{spherical caps}. When the spherical caps are exactly \emph{hemispheres} (pointing in uniformly random directions), the critical threshold $\alpha_c=2$ was determined by combinatorial arguments \cite{cover1965geometrical,MR146858}. Later work showed that as long as the spherical caps have \emph{up to half} the volume of the entire sphere, the limiting free energy and critical threshold coincide with the physics predictions \cite{MR1964377,stojnic2013another}. These results rely crucially on the \textbf{convex} nature of the model. In the regime where the spherical caps have \emph{more than half} the volume of the entire sphere, the problem is no longer convex, and the expected behavior is more complicated \cite{franz2016simplest,franz2017universality,alaoui2020algorithmic}.

For \textbf{Ising} perceptron models, in comparison with the spherical versions, rather less has been rigorously proved. Many of the relevant works were already mentioned at the start of the section. There has been a very fruitful line of works investigating \textbf{symmetric} perceptron models, which are much more tractable because the \textbf{moment method} typically gives fairly sharp results \cite{MR3983947,perkins2021frozen,abbe2021proof,abbe2021binary,gamarnik2022algorithms}. Several of these recent works have investigated algorithmic properties of the solution landscape, inspired in part by conjectures in the physics literature (see e.g.\ \cite{baldassi2016unreasonable}). In particular, for a symmetric interval perceptron model --- more precisely, the model \eqref{e:Z.U} with 
$U(x)=\Ind{x\in[-\kappa,\kappa]}$ and gaussian or Bernoulli disorder $\XI$ --- the sharp threshold and free energy were determined by \cite{perkins2021frozen,abbe2021proof}. The symmetric interval model is also related to the discrepancy minimization problem; see e.g.\ \cite{MR784009,MR3024770,MR3416145,turner2020balancing,MR4398820}.

For non-symmetric Ising perceptron models, it is much more difficult to understand the typical size or behavior of the perceptron solution set $S$. The moment method does not give sharp results. Bounds on the critical window were given by \cite{MR1629627,MR1716771}; and the free energy was computed for small $\alpha$ by \cite{MR1782273,MR3024566,bolthausen2021gardner}. A result of \cite{ding2018capacity} shows that, under an additional condition, the random set $\bS$ is nonempty with \emph{positive} probability throughout the predicted regime; our result allows to improve this statement to high probability. For the half-space model with Bernoulli disorder, \cite{MR4317708} proved a sharp threshold sequence using a characterization of Hatami for low-influence boolean functions \cite{MR2925389}. It would be of interest to see if the methods of \cite{MR4317708} can be extended to more general perceptron models, in ways that do not require more precise estimates on these models \cite{changjixu}. With respect to our current paper, the most closely related previous results are the estimates obtained by Talagrand for the half-space perceptron model \cite{MR1680236,MR3024566}.

\subsection*{Organization} The remainder of this paper is organized as follows:

\begin{itemize}
\item In Section~\ref{s:prelim} we review Talagrand's  
\hyperlink{proof:t.hspace.gaus}{proof of Theorem~\ref{t:hspace.gaus}}, for the half-space perceptron \eqref{e:Z.hspace} with gaussian disorder. We indicate the obstructions to extending this result to more general disorder distributions.

\item In Section~\ref{s:interval} we prove our main ``add one constraint'' estimates, 
Theorems~\ref{t:hspace.general}--\ref{t:interval.general}, for 
the models \eqref{e:Z.hspace} and \eqref{e:Z.hspace} with general disorder.

\item In Section~\ref{s:univ} we give a weak form of the multivariate central limit theorem, and use it to prove an averaged universality statement, Theorem~\ref{t:univ.avg}, for the perceptron models 
\eqref{e:Z.hspace} and \eqref{e:Z.U}.

\item In Section~\ref{s:sharp} we prove Theorems~\ref{t:conc}--\ref{t:univ}.
\end{itemize}

\subsection*{Acknowledgements}  We wish to thank David Belius, Erwin Bolthausen, Ryoki Fukushima, Elchanan Mossel, Joe Neeman, and Changji Xu for many interesting conversations. S.N.\ is supported in part by SNSF grant 176918. N.S.\ is supported in part by NSF CAREER grant DMS-1940092 and NSF-Simons grant DMS-2031883.

\section{Talagrand's results on the half-space perceptron}\label{s:prelim}

In this section we review Talagrand's results on the half-space perceptron \eqref{e:Z.hspace}. We begin with the case of gaussian disorder, and then discuss the possibility of extending to more general distributions of $\XI$: 
\begin{itemize}
\item In \S\ref{ss:talagrand} we review Talagrand's \hyperlink{proof:t.hspace.gaus}{proof of Theorem~\ref{t:hspace.gaus}}.
We highlight two key ingredients in the proof, Lemma~\ref{l:sudakov} (Sudakov minoration) and Lemma~\ref{l:borell.tis} (Borell--TIS).
\item In \S\ref{ss:sup.conc} we discuss extensions of
Lemma~\ref{l:borell.tis} to more general distributions.
\item In \S\ref{ss:sudakov} we discuss extensions of
Lemma~\ref{l:sudakov} to more general distributions.

\end{itemize}
Throughout this paper we use $C$, $\bar{C}$, $C'$, and $C''$ to denote absolute constants. Following common convention, the value of the constant may change from one occurrence to the next, but in a way that does not depend on $N$ or any of the parameters of the model. We use the different labels $C$, $\bar{C}$, $C'$, and $C''$ to avoid ambiguities when different constants interact in the same proof.

\subsection{Intersection of cube and half-space with gaussian disorder}
\label{ss:talagrand}

In this subsection we review some results from \cite{MR1680236,MR3024566}, including the \hyperlink{proof:t.hspace.gaus}{proof of Theorem~\ref{t:hspace.gaus}}. As commented above, the paper \cite{MR1680236} considers the model \eqref{e:Z.hspace} where the sum goes over $x\in\set{-1,+1}^N$, and the $g^a$ are replaced with $\XI^a$ which are i.i.d.\ uniform from $\set{-1,+1}^N$. Meanwhile \cite[Ch.~9]{MR3024566} treats the gaussian case, but presents a simplified argument that does not appear to handle the small $\delta$ regime. We present a more complete summary of the gaussian case below. We begin by recalling two well-known results:

\begin{lem}[Sudakov minoration for gaussian processes]\label{l:sudakov}
Let $(u_i)_{i\le n}$ be a centered gaussian process with $\E[(u_i)^2]=1$ for all $i$, and $\E(u_iu_j)\le 1-\epsilon<1$ for all $i\ne j$. Then
	\[\E\bigg(\max_{i\le n} u_i\bigg)
	\ge \f{(\epsilon\log n)^{1/2}}{2^{1/2}}\,.\]
\begin{proof}See for example \cite[Thm.~13.4]{MR3185193}.
This is discussed further in \S\ref{ss:sudakov} below.
\end{proof}
\end{lem}

\begin{lem}[Borell--TIS inequality]
\label{l:borell.tis}
Let $(u_i)_{i\le n}$ be a centered gaussian process with $\E[(u_i)^2]=1$ for all $i$. Let $u_{\max}\equiv\max_{i\le n} u_i$. Then, for any $s\ge0$, we have
	\[
	\max\bigg\{
	\P\Big( u_{\max} -\E u_{\max} \ge s\Big),
	\P\Big( u_{\max} -\E u_{\max} \le -s\Big)
	\bigg\}
	\le
	\f1{\exp(s^2/2)}\,.
	\]
\begin{proof}
See for example \cite[Thm.~5.8]{MR3185193}. The concentration for $u_{\max}$ is the consequence of a more general concentration result for Lipschitz functionals of gaussian processes; see \cite[Thm.~5.6]{MR3185193}. This is discussed further in \S\ref{ss:sudakov} below.
\end{proof}
\end{lem}

An immediate consequence of the two preceding lemmas is the following:

\begin{cor}\label{c:all.fail}
Let $(u_i)_{i\le n}$ be a centered gaussian process with $\E[(u_i)^2]=1$ for all $i$, and $\E(u_iu_j)\le 1-\epsilon<1$ for all $i\ne j$. 
Let $u_{\max}\equiv\max_{i\le n} u_i$. Then
	\[
	\P\bigg(
	u_{\max} \le \f{(\epsilon\log n)^{1/2}}{2}
	\bigg) \le \f1{n^{\epsilon/50}}\,.
	\]
\begin{proof}
It follows by Lemma~\ref{l:sudakov} combined with Lemma~\ref{l:borell.tis} that
	\[\P\bigg(u_{\max} \le \f{(\epsilon\log n)^{1/2}}{2}\bigg)
	\le \P\bigg( u_{\max}
	-\E u_{\max} \le - \bigg(\f1{2^{1/2}}-\f12\bigg) (\epsilon\log n)^{1/2}\bigg)
	\le \f1{n^{\epsilon/50}}\,,
	\]
as claimed.
\end{proof}
\end{cor}

We next record some basic notations which will be used throughout:

\begin{dfn}
For $t=1-\epsilon$, denote the binary relative entropy function
	\beq\label{e:binary.relent} k_2(t)
	\equiv H\bigg(\f{1+t}{2}\,\bigg|\,\f12\bigg)
	= \f{1+t}{2}\log(1+t)
	+\f{1-t}{2}\log(1-t)
	\equiv \log2-\psi_2(\epsilon)
	,\eeq
so $0\le k_2(t)\le\log2$ with $k_2(0)=0$ and $k_2(-1)=k_2(1)=\log2$.
It will also be useful to recall a simplified bound:
for $0\le p\le 1$ and $0\le tp\le1$ we have
	\beq\label{e:chernoff.simplified}
	H( tp\,|\,p)
	= tp\log t
	+ (1-tp)\log\f{1-tp}{1-p}
	\ge tp\log t+(1-tp)\log(1-tp)
	\ge tp\log \f{t}{e}\,,
	\eeq
where the last step uses that $(1-x)\log(1-x)\ge -x$ for all $x\le1$. (Note however that since the relative entropy is always nonnegative, the bound \eqref{e:chernoff.simplified} is vacuous unless $t\ge e$.)
\end{dfn}

The following estimate says that in any large subset of $\set{-1,+1}^N$, most pairs of elements of that subset are separated in Hamming distance, where the separation guarantee depends on the size of the set:

\begin{lem}[rephrasing of {\cite[Lem.~2.2]{MR1680236}}]
\label{l:inner.prod} Let $S$ be any subset of $\set{-1,+1}^N$ with $|S|\ge\exp(N\delta)$. If $\mu$ denotes the uniform probability measure on $S$, then
	\[
	\mu^{\otimes 2}\bigg( \bigg|\f{(\sigma^{(1)},\sigma^{(2)})}{N}\bigg|
	\ge 1-\epsilon\bigg)
	\le \f{2}{\exp(N\delta/2)}\,,
	\]
for any $\epsilon\in(0,1)$ with $\psi_2(\epsilon)\le\delta/2$.\footnote{Note this result is similar to but stronger than \cite[Lem.~9.2.1]{MR3024566}.}
\begin{proof} Abbreviate $t\equiv 1-\epsilon$. Conditioning on $\sigma^{(1)}$, we have
	\begin{align*}
	&\mu\bigg(\bigg\{\sigma^{(2)} : 
	\bigg|\f{(\sigma^{(1)},\sigma^{(2)})}{N}\bigg| \ge t
		\bigg\}\bigg)
	= \f1{|S|}
		\bigg|
		\bigg\{\sigma^{(2)}\in\set{-1,+1}^N :  \bigg|\f{(\sigma^{(1)},\sigma^{(2)})}{N}\bigg| \ge t
		\bigg\}\bigg| \\
	&\le \f{2^N}{\exp(N\delta)} \P
	\bigg( \bigg|
	\textup{Bin}\bigg(N,\f12\bigg) -\f{N}{2}\bigg| \ge \f{Nt}{2}
	\bigg)
	\le \f{2^{N+1} \exp(-Nk_2(t))}{\exp(N\delta)}
	\le \f{2\exp(N\psi_2(\epsilon))}{\exp(N\delta)}
	\,.
	\end{align*}
The claim follows.
\end{proof}
\end{lem}

The next result is the main ingredient in Talagrand's proof of Theorem~\ref{t:hspace.gaus}:

\begin{ppn}[{adaptation of \cite[Propn.~2.1]{MR1680236}}]
\label{p:largevals.gaus}
Let $S$ be any subset of $\set{-1,+1}^N$ with $|S|\ge\exp(N\delta)$ where $\delta$ is a positive constant. Let $\epsilon\in(0,1)$ with $\psi_2(\epsilon)\le\delta/2$. Then
	\[
	\P\bigg(
	\f1{|S|}
	\bigg|
	\bigg\{\sigma\in S : \f{(g,\sigma)}{N^{1/2}} 
	\ge s\epsilon^{1/2} 
	\bigg\}\bigg|
	\le \f1{\exp(C' s^2)}
	\bigg) 
	\le \f{1}{\exp(s^2\epsilon/C')}\,,
	\]
provided $C' \le s \le (N\delta)^{1/2}/C'$, where $C'$ is an absolute constant.

\begin{proof}
Let $\mu$ denote the uniform probability measure on $S$. The basic idea of the proof is the following. By Lemma~\ref{l:inner.prod}, $\mu^{\otimes n}$ gives large weight to well-separated $n$-tuples of configurations. By Corollary~\ref{c:all.fail},
 given a well-separated $n$-tuple, it is very unlikely that all $n$ configurations violate the new constraint. This implies the desired bound. The details are as follows.\medskip

\noindent\textit{Step 1.} First we show that $\mu^{\otimes n}$ gives large weight to well-separated $n$-tuples of configurations. Let $t=1-\epsilon>0$ with $\psi_2(\epsilon)\le\delta/2$ as in Lemma~\ref{l:inner.prod} above. For any positive integer $n$, let $C_n$ denote the subset of configurations in $(\set{-1,+1}^N)^n$ that are pairwise well-separated:
	\[
	C_n
	= \bigg\{(\sigma^{(1)},\ldots,\sigma^{(n)}) \in(\set{-1,+1}^N)^n
		: \f{(\sigma^{(i)},\sigma^{(j)})}{N} \le t \ \forall i\ne j\bigg\}\,.
	\]
It follows from Lemma~\ref{l:inner.prod} that $\mu^{\otimes n}(C_n)$ must be large: indeed, a union bound over all pairs $i<j$ gives
	\beq\label{e:well.sep.lbd}
	\mu^{\otimes n}(C_n)
	\ge 1-\binom{n}{2}\mu^{\otimes 2}\bigg(\bigg\{
		(\sigma^{(1)},\sigma^{(2)})
		:  \bigg|\f{(\sigma^{(1)},\sigma^{(2)})}{N}\bigg| \ge t
		\bigg\}\bigg)
	\ge1-\binom{n}{2} \f{2}{\exp(N\delta/2)}
	\ge\f34\,,
	\eeq
where the last inequality holds provided that
	\beq\label{e:required.ubd.n}
	1\le n \le \f{\exp(N\delta/4)}{2}\,.
	\eeq

\noindent\textit{Step 2.} Now let $H_g$ denote the set of all $\sigma$ such that $(g,\sigma)/N^{1/2}\ge (\epsilon\log n)^{1/2}/2$; we now lower bound the probability that $\mu(H_g)$ is too small. Let $D_{n,g}$ denote the subset of configurations in $(\set{-1,+1}^N)^n$ where \emph{all} $n$ points lie \emph{outside} $H_g$:
	\[
	D_{n,g} = \bigg\{(\sigma^{(1)},\ldots,\sigma^{(n)})
		\in(\set{-1,+1}^N)^n
		: \f{(g,\sigma^{(i)})}{N^{1/2}} < \f{(\epsilon\log n)^{1/2}}{2} 
		 \ \forall 1\le i\le n
		 \bigg\}
		= \Big(\set{-1,+1}^N \setminus H_g\Big)^n
		\,.
	\]
On the event that $\mu(H_g) \le 1/(4n)$, 
the complement of $D_{n,g}$ has measure at most $1/4$ under $\mu^{\otimes n}$, so
	\beq\label{e:all.configs.violate.constraint}
	\f34
	\stackrel{\eqref{e:well.sep.lbd}}{\le} 
	\mu^{\otimes n}(C_n)
	\le
	\mu^{\otimes n}(C_n \cap D_{n,g}) + \f14\,.\eeq
Rearranging the above gives $\mu^{\otimes n}(C_n \cap D_{n,g}) \ge 1/2$ on the event $\mu(H_g) \le 1/(4n)$. It follows that
	\[
	\P\bigg(  \mu(H_g) \le \f{1}{4n} \bigg)
	\stackrel{\eqref{e:all.configs.violate.constraint}}{\le}
	\P\bigg( \mu^{\otimes n}(C_n \cap D_{n,g}) \ge \f12 \bigg)
	\le 2\E\Big[\mu^{\otimes n}(C_n \cap D_{n,g})\Big]
	\le \f{2}{n^{\epsilon/50}}\,,
	\]
where the intermediate step is by Markov's inequality, and the last step is by Corollary~\ref{c:all.fail}. Making a change of variables $s^2 = (\log n)/4$ gives the claim.
\end{proof}\end{ppn}

\begin{proof}[\hypertarget{proof:t.hspace.gaus}{Proof of Theorem~\ref{t:hspace.gaus}}] If $C' \le s \le (N\delta)^{1/2}/C'$ and $s\epsilon^{1/2}\ge \kappa$, then Proposition~\ref{p:largevals.gaus} immediately implies
	\[
	\P\bigg(
	\f1{|S|}
	\bigg|
	\bigg\{\sigma\in S : \f{(g,\sigma)}{N^{1/2}} \ge  \kappa
	\bigg\}\bigg|
	\le \f1{\exp(C' s^2)}
	\bigg) 
	\le \f{1}{\exp(s^2\epsilon/C')}\,.
	\]
The claim follows. Note that if $\kappa\le0$ then the condition $s\epsilon^{1/2}\ge\kappa$ holds trivially for all $s\ge0$, so the $C_\delta$ in the theorem statement can be taken independently of $\kappa$. By contrast, if $\kappa>0$ then $C_\delta$ depends on $\kappa$ as well.
\end{proof}

For comparison, the next proposition is a variant of Proposition~\ref{p:largevals.gaus} which was also proved by Talagrand but with different methods. It gives similar although not exactly comparable results. In this paper we will only make use of  Proposition~\ref{p:largevals.gaus}, so Proposition~\ref{p:largevals.gaus.TAL} can be skipped by the reader.

\begin{ppn}[{adaptation of \cite[Thm.~8.2.4]{MR3024566}}]
\label{p:largevals.gaus.TAL}
Let $(u_i)_{i\le n}$ be a centered gaussian process with $\E[(u_i)^2]=1$ for all $i$, and 
$\E[u_i u_j]\le 1-\epsilon$ for all $i\ne j$. Then
	\[
	\P\bigg(
	\f1n\Big|
	\Big\{i\le n: u_i \ge s\epsilon^{1/2}\Big\}
	\Big| <
	\f1{\exp(\bar{C}s^2/\epsilon)}
	\bigg) 
	\le \f1{\exp(s^2/2)}
	\]
provided $\bar{C} \le s \le (\log n)^{1/2}/\bar{C}$, where $\bar{C}$ is an absolute constant.

\begin{proof}
Fix $s>0$ and consider the functions
	\[
	F(u) 
	\equiv
	\log\bigg( \sum_{i\le n} \exp(su_i)\bigg)\,,\quad
	F_2(u)
	\equiv
	\log\bigg( \sum_{i\le n} \exp(2su_i)\bigg)
	\]
It follows from \cite[Propn.~8.2.2]{MR3024566} (based on \cite[Lem.~8.2.1]{MR3024566})  that if $(u_i)_{i\le n}$ and $(v_i)_{i\le n}$ are two centered gaussian processes with
$\E[(u_i)^2] \ge \E[(v_i)^2]$ for all $i$, and
$\E(u_i u_j)\le\E(v_i v_j)$ for all $i\ne j$, then
$\E F(u) \ge \E F(v)$. In the current setting we can take $v_i = (1-\epsilon)^{1/2}z + \epsilon^{1/2} z_i$ where $z,z_i$ are i.i.d.\ standard gaussian random variables. It follows by combining with \cite[Lem.~8.2.3]{MR3024566} that
	\[
	\E F(u)
	\ge \E F(v)
	=\log\bigg( \sum_{i\le n} \exp(s\epsilon^{1/2}z_i)\bigg)
	\ge \log n + \f{s^2\epsilon}{5}\,,
	\]
where the last inequality holds provided $C \le s \epsilon^{1/2} \le (\log n)^{1/2}/C$, where $C$ is an absolute constant. Moreover, gaussian concentration gives
	\[
	\P\Big( F(u) \le \E F(u) - t\Big)
	\le \exp\bigg(-\f{t^2}{4s^2}\bigg)\,,
	\]
and combining the last two bounds implies
	\[
	\P\bigg( 
	F(u) \ge \log n + \f{s^2\epsilon}{10}\bigg)
	\le 
	1-\exp\bigg(-\f{s^2\epsilon}{400} \bigg)\,.
	\]
On the other hand, it follows using Markov's inequality that
	\[
	\P\bigg( F_2(u) 
	\ge \log n + 3s^2\bigg)
	\le \f{\E\exp F_2(u)}{n\exp(3s^2)}
	= \f{\exp(2s^2)}{\exp(3s^2)}
	= \f1{\exp(s^2)}\,.
	\]
It follows by combining the last two bounds that the event
	\[
	E_\star
	\equiv \bigg\{ 
	F(u) \ge \log n + \f{s^2\epsilon}{10},
	F_2(u) 
	\le\log n + 3s^2
	\bigg\}
	\]
has probability $\P(E_\star)\ge 1-2\exp(-s^2\epsilon/400)$. Let $\mathbf{P}_n$ be the uniform probability measure on $\set{1,\ldots,n}$, and on the probability space $([n],\mathbf{P}_n)$
consider the random variable $X : i \mapsto  X(i) = \exp(su_i)$. Then
$\mathbf{E}_n X = n^{-1}\exp F(u)$ and $\mathbf{E}_n(X^2) =n^{-1} \exp F_2(u)$, so the Paley--Zygmund inequality implies
	\beq\label{e:paleyzyg}
	\mathbf{P}_n\bigg(
	X\ge \f{\mathbf{E}_n X}{2}\bigg)
	\ge \f{(\mathbf{E}_n X)^2}{4\mathbf{E}_n(X^2)}
	=\f{ (n^{-1}\exp F(u))^2}{4 n^{-1} \exp F_2(u)}
	\ge \f{\exp(s^2\epsilon/5)}{4\exp(3s^2)}
	\ge \f1{\exp(3s^2)}\,,
	\eeq
where the second-to-last inequality holds on the event $E_\star$, and the last inequality holds using
$s\epsilon^{1/2}\ge C$. Moreover,
	\[
	\f{\mathbf{E}_n X}{2}
	= \f{\exp F(u)}{2n}
	\ge \f{\exp(s^2\epsilon/10)}{2}
	\ge \exp\bigg( \f{s^2\epsilon}{20}\bigg)\,,
	\]
where the last inequality again uses
$s\epsilon^{1/2}\ge C$. Thus,  $X(i) = \exp(su_i) \ge \mathbf{E}_n X/2$ implies $u_i \ge s\epsilon/20$, so
	\[
	\mathbf{P}_n\bigg(
	X\ge \f{\mathbf{E}_n X}{2}\bigg) \le
	\f1n\bigg|\bigg\{ i\le n : 
	u_i \ge \f{s\epsilon}{20}\bigg\}\bigg|\,.
	\]
Since the bound \eqref{e:paleyzyg} holds on the event $E_\star$, it follows that 
	\begin{align*}
	\P\bigg(\f1n\bigg|\bigg\{ i\le n : 
	u_i \ge \f{s\epsilon}{20}
	\bigg\}\bigg| < \f1{\exp(3s^2)}
	\bigg)
	&\le \P\bigg(
	\mathbf{P}_n\bigg(
	X\ge \f{\mathbf{E}_n X}{2}\bigg)
	< \f1{\exp(3s^2)}
	\bigg) \\
	&\le 1-\P(E_\star) 
	\le 2\exp\bigg(-\f{s^2\epsilon}{400}\bigg)
	\le \exp\bigg(-\f{s^2\epsilon}{800}\bigg)
	\,.\end{align*}
Taking $s'=s\epsilon^{1/2}/20$ gives, for $C \le 20 s'  \le (\log n)^{1/2}/C$,
	\[
	\P\bigg(\f1n\bigg|\bigg\{ i\le n : 
	u_i \ge s'\epsilon^{1/2}
	\bigg\}\bigg| 
	< \f1{\exp(1200(s')^2/\epsilon)}
	\bigg)
	\le \exp\bigg(-\f{(s')^2}{2}\bigg)\,.
	\]
The claimed bound follows by setting $\bar{C}=\max\set{1200,20C}$.
\end{proof}
\end{ppn}

\subsection{Concentration of supremum of canonical processes}
\label{ss:sup.conc}

In the remainder of this section we discuss the possibility of extending the approach of Theorem~\ref{t:hspace.gaus} to more general distributions of the disorder $\XI$. The main conclusion of this discussion is that the proof of Theorem~\ref{t:hspace.gaus} can likely be extended to a limited class of distributions using results from the existing literature. By contrast, our main results follow a different (and perhaps simpler) approach, and apply to the wider class of all subgaussian distributions. However,  we obtain somewhat weaker bounds than in Theorem~\ref{t:hspace.gaus} (e.g., 
Theorem~\ref{t:hspace.gaus} applies for $w\le N/C_\delta$ while
Theorem~\ref{t:hspace.general} applies for $w\le N^{1/2}/C_\delta$).

Theorem~\ref{t:hspace.gaus} can be obtained as a consequence of either Proposition~\ref{p:largevals.gaus} or Proposition~\ref{p:largevals.gaus.TAL}. Proposition~\ref{p:largevals.gaus.TAL} relies on a gaussian interpolation bound and appears difficult to extend more generally, so we 
turn to discussing the possibility of extending  Proposition~\ref{p:largevals.gaus}. The proposition relies on two essential ingredients:
\begin{enumerate}[(i)]
\item The Sudakov minoration (Lemma~\ref{l:sudakov}), lower bounding
the expected supremum of a gaussian process;
\item The Borell--TIS inequality (Lemma~\ref{l:borell.tis}), giving concentration of the supremum of a gaussian process.
\end{enumerate}
We discuss each of these ingredients separately below, for a stochastic process
	\beq\label{e:canonical}
	\bigg(\f{(\XI,\sigma)}{N^{1/2}} : \sigma\in S\bigg)\,,
	\eeq
where $S$ is a subset of $\R^N$. The rescaling of $S$ by $N^{1/2}$ is not essential; we chose it to maintain consistency with the scaling in the rest of this paper, where $S$ is generally a subset of $\set{-1,+1}^N$. In the literature, the 
stochastic process \eqref{e:canonical} is sometimes termed the ``canonical process'' for $\XI$ indexed by $S/N^{1/2}$. In this subsection we discuss concentration of the supremum of the process (point (ii) above); in \S\ref{ss:sudakov} we discuss lower bounding the expected supremum (point (i) above).
 
For $S\subseteq\R^N$, let us consider the supremum of the canonical process \eqref{e:canonical} as a function of the disorder $\XI$:
	\beq\label{e:sup.canonical}
	f(\XI)
	\equiv f\bigg(\XI; \f{S}{N^{1/2}}\bigg)
	\equiv \sup\bigg\{ \f{(\XI,\sigma)}{N^{1/2}} : \sigma\in S\bigg\}\,.
	\eeq
If $S$ is any subset of $N^{1/2} B_2$ where $B_2$ is the euclidean unit ball in $\R^N$, then $f:\R^N\to\R$ will be Lipschitz with respect to the euclidean norm on $\R^N$:
	\beq\label{e:sup.Lip}
	f(\XI')
	= \sup\bigg\{ \f{(\XI,\sigma)}{N^{1/2}}
		+\f{(\XI'-\XI,\sigma)}{N^{1/2}} : \sigma\in S\bigg\}
	\le f(\XI) + \|\XI'-\XI\|_2
	\eeq
for any $\XI,\XI'\in\R^N$. The Borell--TIS inequality gives concentration for any Lipschitz functional of a gaussian random variable, hence Lemma~\ref{l:borell.tis} (since the gaussian process $(u_i)_{i\le n}$ can be expressed as $((g,v_i))_{i\le n}$, where $g$ is a standard gaussian vector and $(v_i)_{i\le n}$ is a collection of unit vectors). We know of similar bounds in two other settings:
\begin{itemize}
\item Suppose $\XI_i$ has density $\exp(-V(x))\,dx$ 
where $V$ is \textbf{uniformly convex} (i.e., the second derivative of $V$ is uniformly bounded below by a positive constant). Then it follows from \cite[Propn.~2.18]{MR1849347} (see also \cite[Thm.~5.2.15]{MR3837109}) that \eqref{e:sup.Lip} satisfies a similar concentration bound
as Lemma~\ref{l:borell.tis}, for any $S\subseteq N^{1/2}B_2$.
\item If the $\XI_i$ are almost surely \textbf{bounded} random variables, then
it can be deduced via the method of bounded differences (see e.g.\ \cite[Thm.~6.2]{MR3185193}) that
\eqref{e:sup.Lip} satisfies a similar concentration bound
as Lemma~\ref{l:borell.tis}, again for any $S\subseteq N^{1/2}B_2$.
\end{itemize}
For subgaussian distributions, we do not have such strong guarantees, but we have the following lemma which gives a weaker result under the further assumption $S\subseteq\set{-1,+1}^N$ (rather than $S\subseteq N^{1/2}B_2$):

\begin{lem}\label{l:conc.sup.subgaussian}
Suppose $\XI$ is a random vector in $\R^N$ satisfying Assumption~\ref{a:subgaus}. If $S\subseteq \set{-1,+1}^N$, then
	\[
	\P\bigg( \bigg| f\bigg(\XI;\f{S}{N^{1/2}}\bigg)
		 - \E f\bigg(\XI;\f{S}{N^{1/2}}\bigg)
		 \bigg| \ge u \bigg)
	\le \exp\bigg( - \min\bigg\{ \f{u^2}{4 A^2},
	\f{N^{1/2}u}{2A} \bigg\}\bigg)
	\]
for $f$ as defined by \eqref{e:sup.canonical}, and $A$ a constant depending only on $\NU$.

\begin{proof}
Let $\HH_i$ denote the $\sigma$-algebra generated by the first $i$ coordinates $(\XI_1,\ldots,\XI_i)$, and decompose
	\[
	f(\XI)-\E f(\XI)
	= \sum_{i=1}^N (X_i-X_{i-1})
	\]
where $X_i\equiv \E(f(\XI)\,|\,\HH_i)$. 
Let $\E_i$ be expectation over $\XI_i$ only,
and let $\E^i$ be expectation over all the $(\XI_j)_{j\ge i}$, so that we have $X_i = \E^{i+1} f(\XI)$ while $X_{i-1}=\E^i f(\XI) = \E_i \E^{i+1} f(\XI)$. We then have, using Jensen's inequality,
	\begin{align*}
	E_i(\lambda) &\equiv
	\E\bigg[ \exp\Big(\lambda|X_i-X_{i-1}|\Big)
	\,\bigg|\,\HH_{i-1}\bigg]
	=\E_i \E^{i+1}
	\exp\bigg( \lambda\Big| \E^{i+1}\Big(f(\XI)-\E_if(\XI)\Big)\Big|\bigg)
	\\
	&
	\le \E_i \E^{i+1}
	\exp \bigg( \lambda\Big| f(\XI)-\E_if(\XI)\Big|\bigg)\,,
	\end{align*}
Let $\XI^{(i)}$ denote the vector that results from replacing the $i$-th coordinate of $\XI$ with an independent copy $\zeta_i$, and let 
$\E^{(i)}$ denote expectation over $\zeta_i$ only. Then applying Jensen's inequality again gives
	\[
	E_i(\lambda)
	\le
	\E_i \E^{i+1}
	\exp \bigg( \lambda\Big|
	\E^{(i)}\Big(f(\XI)-f(\XI^{(i)})\Big)
	 \Big|\bigg)
	\le \E_i \E^{i+1}\E^{(i)}
	\exp \bigg( \lambda\Big| f(\XI)-f(\XI^{(i)})\Big|\bigg)\,.\]
We then note, since we assume $S\subseteq\set{-1,+1}^N$, 
the bound \eqref{e:sup.Lip} can be refined to
	\[
	f(\XI^{(i)})
	\le f(\XI)
	+ \max\bigg\{\f{(\XI'-\XI,\sigma)}{N^{1/2}} : \sigma\in \set{-1,+1}^N\bigg\}
	\le f(\XI) + \f{|\zeta_i-\XI_i|}{N^{1/2}}\,.
	\]
Since for any $x\in\R$ we have $\exp(|x|)\le \exp(x) + \exp(-x)$, we can use Assumption~\ref{a:subgaus} to bound
	\[
	E_i(\lambda)
	\le \E\exp\bigg(\f{\lambda|\zeta_i-\XI_i|}{N^{1/2}}\bigg)
	\le 2\exp\bigg( \f{\lambda^2\NU}{N}\bigg)\,.
	\]
In particular, taking $\lambda=N^{1/2}$ gives $E_i(N^{1/2}) \le 2\exp(\NU)$. Another application of Jensen's inequality gives that for $A$ large enough  (depending on $\NU$ only), we have 
	\[
	E_i\bigg( \f{N^{1/2}}{A}\bigg) 
	\le \Big( E_i(N^{1/2}) \Big)^{1/A} 
	\le \Big( 2\exp(\NU)\Big)^{1/A} \le 2\,.
	\]
It then follows by the martingale Bernstein inequality (see \cite[Thm.~A.6.1]{MR3024566}) that
	\[
	\P\bigg(\Big| f(\XI)-\E f(\XI)\Big|\ge u\bigg)
	\le
	\exp\bigg( - \min\bigg\{ \f{u^2}{4 A^2},
	\f{N^{1/2}u}{2A} \bigg\}\bigg)\]
for all $u\ge0$, as claimed.
\end{proof}
\end{lem}

We note that, for the purposes of Theorem~\ref{t:hspace.gaus},
Lemma~\ref{l:conc.sup.subgaussian} can be substituted for Lemma~\ref{l:borell.tis} with similar results. Thus, in our setting, concentration of the supremum of \eqref{e:canonical} does not appear to be a major issue. Rather, as we next discuss, lowering bounding the expected supremum appears to be the main obstacle to extending the argument of Theorem~\ref{t:hspace.gaus} to more general distributions.\medskip

\subsection{Sudakov minoration lower bound}
\label{ss:sudakov}

We now discuss lower bounding the expected supremum
	\beq\label{e:canonical.sup.mm}
	\E f\bigg( \XI,\f{S}{N^{1/2}}\bigg)
	= \E\bigg[
	\sup\bigg\{ \f{(\XI,\sigma)}{N^{1/2}} : \sigma\in S
	\bigg\}\bigg]\,,\eeq
for $f$ as in \eqref{e:sup.canonical} and $S\subseteq\R^N$. In the gaussian case, the desired lower bound is given by Lemma~\ref{l:sudakov}. In the \textbf{Bernoulli} case a comparable lower bound is given by \cite[Propn.~2.2]{MR1207231}. Minoration bounds for more general processes do exist (e.g.\ \cite{MR1207231,MR1269606,MR1475551,MR2133757}), but the dependence on the law of $\XI_i$ is somewhat complicated. The most relevant results come from \cite{MR1475551}
(see also \cite[Ch.~5]{MR2133757}), for the situation that the $\XI_i$ are symmetric random variables such that the tail probability function
	\beq\label{e:tail.prob.fn}
	\phi_\XI(x)
	\equiv \log \f1{\P(|\XI_i|\ge x)}\eeq
is \textbf{convex}. For simplicity, we will restrict most of our discussion to the case that the distribution of $\XI_i$ has density
	\beq\label{e:stable.law}
	\f{\exp(-|x|^\alpha)}{z_\alpha}\,dx
	\eeq
with respect to Lebesgue measure, where $1\le\alpha<\infty$ and $z_\alpha$ denotes the normalizing constant. This is a special case of the setting of \cite{MR1475551}, and was treated earlier by \cite{MR1207231,MR1269606}.

If $\XI_i$ follows the distribution \eqref{e:stable.law}
with $1\le\alpha\le2$, or more generally if $\phi_\XI(x)\lesssim x^2$ for $x\ge1$ (meaning that $\XI_i$ has \textbf{heavier tails} than the gaussian distribution), then it follows from \cite[Thm.~1]{MR1475551} (see also \cite[Thm.~1.2]{MR1269606} and the discussion around \cite[Thm.~5.2.6]{MR2133757}) that a similar lower bound as Lemma~\ref{l:sudakov} holds.
In more detail, for an $N$-dimensional gaussian vector $g$ and for
any $X\subseteq\R^N$, it was proved by \cite{MR906527} that
	\[\E\bigg[\sup\Big\{ (g,x) : x\in X \Big\}
	\bigg]
	\asymp
	\gamma_2(X,\|\cdot\|_2)
	=\adjustlimits \inf_{(X_n)}
	\sup_{t\in X}
	\sum_{n\ge0} 2^{n/2} 
	\min\Big\{\|t-x\|_2:x\in X_n\Big\}
	\,,
	\]
where $X_n$ is a subset of $X$ of cardinality at most $2^{2^n}$
($\gamma_2$ is a standard notation in the topic of majorizing measures or generic chaining). If the $\XI_i$ are symmetric random variables with $\phi_\XI(x) \lesssim x^2$ for $x\ge1$, then we have
	\[ \E\bigg[
	\sup\Big\{ (\XI,x) : x\in X\Big\}
	\bigg]
	\gtrsim \gamma_2(X,\|\cdot\|_2)
	\asymp
	\E\bigg[
	\sup\Big\{ (g,x) : x\in X\Big\}
	\bigg]
	\] 
for any $X\subseteq\R^N$ (cf.\ and \cite[Thm.~2.1.1 and display~(5.44)]{MR2133757}). If $X$ is a subset of the unit sphere in $\R^N$ such that $(x,x')\le 1-\epsilon$ for all $x\ne x'$ in $X$, then the right-hand side of the above is lower bounded by Lemma~\ref{l:sudakov}.\footnote{The resulting lower bound on the left-hand side may not be tight when $\alpha<2$, although this in itself does not appear to be an issue.}

If $\XI_i$ follows the distribution \eqref{e:stable.law}
with $\alpha\ge2$  (so that $\XI_i$ has \textbf{lighter tails} than the gaussian distribution), then characterizing the expected supremum is in general more complicated; see \cite[Thm.~1.3]{MR1269606}. A lower bound that is easier to work with is given by 
\cite[Thm.~3.1]{MR1269606}
(see also \cite[Thm.~1]{MR1475551}) which in fact applies for all $\alpha>1$: writing $\mm\equiv\mm(S/N^{1/2})$ for the quantity in 
\eqref{e:canonical.sup.mm}, and writing 
$\beta$ for the dual (H\"older conjugate) of $\alpha$, we have
	\beq\label{e:sudakov.stable}
	\mm
	\gtrsim\log n\bigg( \f{S}{N^{1/2}},U_\beta(\mm)\bigg)\,,
	\eeq
where $n(X,U)$ denotes the number of translates of $U$ needed to cover $X$, and
	\[
	U_\beta(u)
	\equiv \bigg\{ x : \sum_{i\ge1}\eta_\beta(x_i)
	\le u \bigg\}\,,\quad
	\eta_\beta(x)\equiv\begin{cases} x^2 & \textup{for $|x|\le1$,}\\
		|x|^\beta & \textup{for $|x|\ge1$.}
		\end{cases}
	\]
Adapting the argument following \cite[Thm.~3.2]{MR1269606} yields the following bound, which (like Lemma~\ref{l:conc.sup.subgaussian}) relies crucially on the restriction $S\subseteq\set{-1,+1}^N$:

\begin{cor}\label{c:sudakov.stable}
Suppose $S\subseteq\set{-1,+1}^N$,  and assume that for all $\sigma\ne\tau$ in $S$ we have
$|(\sigma,\tau)|/N \le 1-\epsilon$. Suppose $\XI$ is an $N$-dimensional random vector with i.i.d.\ entries $\XI_i$ distributed according to \eqref{e:stable.law} with $\alpha>1$. Then
	\[ \E\bigg[
	\sup\bigg\{ \f{(\XI,\sigma)}{N^{1/2}} : \sigma\in S
	\bigg\}
	\bigg]
	\gtrsim \epsilon (\log |S|)^{1/2}\,.
	\]

\begin{proof}
As above, we abbreviate $\mm=\mm(S/N^{1/2})$ for the quantity of interest.  Let $\theta$ denote a positive constant.
Applying \eqref{e:sudakov.stable} (from \cite[Thm.~3.1]{MR1269606}) with $\theta S$ in place of $S$ gives
	\[
	\theta \mm\gtrsim
	\log n\bigg(\f{\theta S}{N^{1/2}},
	U_\beta(\theta \mm)
	\bigg)\,.
	\]
Suppose two points $x,y\in[-1/2,+1/2]^N$ are covered by the same translate of $U_\beta(\theta \mm)$. Then there must exist $z\in[-1/2,+1/2]^N$ such that both $x$ and $y$ belong to the translate of $U_\beta(\theta \mm)$ centered at $z$: that is, 
	\[\max\bigg\{
	\sum_{i\le N} \eta_\beta(x_i-z_i),
	\sum_{i\le N} \eta_\beta(y_i-z_i)
	\bigg\}
	\le \theta\mm\,.
	\]
Since $|x_i-z_i|\le1$ and $|y_i-z_i|\le1$ for all $i\le N$, the definition of $\eta_\beta$ implies that both $x$ and $y$ must lie within \emph{euclidean} distance $(\theta\mm)^{1/2}$ of $z$, and so $x$ and $y$ must also be covered by the same translate of $(\theta\mm)^{1/2} B_2$, where $B_2$ denotes the unit ball in euclidean norm. Thus, as long as $\theta/N^{1/2} \le 1/2$, we have
	\[
	\theta\mm\gtrsim
	\log n\bigg(\f{\theta S}{N^{1/2}},
	(\theta\mm)^{1/2} B_2
	\bigg)
	=\log n\bigg( \f{S}{N^{1/2}}, 
	\bigg(\f{\mm}{\theta}\bigg)^{1/2} B_2 \bigg)\,.
	\]
If $\mm/\epsilon \le N^{1/2}/2$, then setting $\theta=\mm/\epsilon$ gives
	\[
	\f{\mm^2}{\epsilon} \ge  \log n\bigg(
	\f{S}{N^{1/2}}, \epsilon^{1/2} B_2\bigg)
	= \log |S|\,,
	\]
since the separation assumption in the statement of the corollary
implies that each translate of $\epsilon^{1/2} B_2$ can cover at most one element of $S$. Otherwise we must have $\mm \ge \epsilon N^{1/2}/2 \gtrsim \epsilon (\log |S|)^{1/2}$, since $|S|\le 2^N$. Combining these two cases gives the claimed lower bound.
\end{proof}
\end{cor}

For simplicity we have stated Corollary~\ref{c:sudakov.stable} for the distribution \eqref{e:stable.law} with $\alpha>1$, so that the bound can be obtained by applying \cite[Thm.~3.1]{MR1269606}. For the more general setting where the function $\phi_\XI(x)$ of \eqref{e:tail.prob.fn} is \textbf{convex}, a similar lower bound as Corollary~\ref{c:sudakov.stable}  can be obtained by applying \cite[Thm.~1]{MR1475551}. In summary, we have similar lower bounds for the canonical process \eqref{e:canonical} in the cases
\begin{itemize}
\item The $\XI_i$ are gaussian (Lemma~\ref{l:sudakov});
\item The $\XI_i$ are  Bernoulli (\cite[Propn.~2.2]{MR1207231});
\item The $\XI_i$ are symmetric random variables such that the function
 $\phi_\XI(x)$ is convex  (\cite[Thm.~3.1]{MR1269606},
\cite[Thm.~1]{MR1475551}, and
Corollary~\ref{c:sudakov.stable}).
\end{itemize}
We do not know of a similar lower bound for the more general class of subgaussian vectors.\medskip

\noindent\textbf{Summary of conclusions of \S\ref{ss:sup.conc}--\ref{ss:sudakov}.} The conclusion of the last two subsections is that,
 using only existing results in the literature, it may be possible to extend Theorem~\ref{t:hspace.gaus} to cover the case that the $\XI_i$ are symmetric subgaussian random variables such that the function $\phi_\XI(x)$ of \eqref{e:tail.prob.fn} is convex. One can adapt the argument of Theorem~\ref{t:hspace.gaus} with 
Corollary~\ref{c:sudakov.stable} in place of Lemma~\ref{l:sudakov}, and
Lemma~\ref{l:conc.sup.subgaussian} in place of Lemma~\ref{l:borell.tis}.
To cover all subgaussian distributions, the main difficulty appears to be in obtaining a lower bound comparable to that of Lemma~\ref{l:conc.sup.subgaussian}. 

In the remainder of this paper we prove versions 
of Theorem~\ref{t:hspace.gaus} for general subgaussian disorder, using as input \textbf{only} the gaussian results Lemma~\ref{l:sudakov} and Lemma~\ref{l:borell.tis}. This results in weaker (and likely suboptimal) bounds. On the other hand, our argument is relatively simple, as it bypasses the more difficult Sudakov minoration bounds proved by \cite{MR1269606,MR1475551}.

\section{Adding a single constraint in perceptron models}
\label{s:interval}

\noindent In this section we prove our main estimates for the perceptron model, Theorems~\ref{t:hspace.general}--\ref{t:interval.general}.
\begin{itemize}
\item In \S\ref{ss:sep.blocks} we prove a preliminary estimate (Proposition~\ref{p:sep.on.blocks.subsample}) 
towards our main results, which says that large subsets of $\set{-1,+1}^N$ must contain many pairs of configurations that are well-separated within blocks of the coordinates $[N]$.
\item In \S\ref{ss:threshold.general} we prove Proposition~\ref{p:largevals.general}, which is a version of Proposition~\ref{p:largevals.gaus} for general disorder. This leads to the \hyperlink{proof:t.hspace.general}{proof of Theorem~\ref{t:hspace.general}}.

\item In \S\ref{ss:addone.general} we apply Proposition~\ref{p:largevals.general} to give the \hyperlink{proof:t.interval.general}{proof of Theorem~\ref{t:interval.general}}. As we explain at the end of that subsection, a similar argument (using Proposition~\ref{p:largevals.gaus} in place of Proposition~\ref{p:largevals.general}) yields the \hyperlink{proof:t.interval.gaus}{proof of Theorem~\ref{t:interval.gaus}}.
\end{itemize}

\subsection{Extraction of subsets with separation on blocks}
\label{ss:sep.blocks}

Recall that we decomposed $[N]$ into blocks $I_1,\ldots,I_L$, each of size $K=N/L$.
Given $S\subseteq\set{-1,+1}^N$ of size $|S|\ge\exp(8N\delta)$, our first goal is to extract $\Omega\subseteq S$ large, such that configurations in $\Omega$ have ``good separation'' on some of the blocks $I_j$. Formally:
\begin{dfn}\label{d:sep}
For any subset of coordinates $I\subseteq[N]$ of size $|I|=K$,
we say that two configurations $\sigma,\tau\in\set{-1,+1}^N$
are \textbf{$\epsilon$-separated on $I$} if
	\[ \f{|(\sigma_I,\tau_I)|}{K} \le 1-\epsilon\,.
	\]
Equivalently, the euclidean distance between $\sigma_I$ and $\tau_I$ 
is at least $(2K\epsilon)^{1/2}$.
\end{dfn}

Note that in an extreme case we could have for instance
	\[
	S = \bigg\{\sigma\in\set{-1,+1}^N
	: \sigma_i = 1 \textup{ for all }
	1\le i \le N\bigg( 1-\f{8\delta}{\log2}\bigg)\bigg\}\,,
	\]
so in this case $|S|\ge\exp(8N\delta)$, but clearly any pair of elements in $S$ will not be well-separated on most blocks $I_j$. Thus it will only be possible to guarantee good separation on a small fraction $\gamma$ (depending on $\delta$) of the blocks $I_j$. The main result of this subsection is the following:

\begin{ppn}\label{p:sep.on.blocks.subsample}
Let $S$ be any subset of $\set{-1,+1}^N$ with $|S|\ge\exp(8N\delta)$ where $\delta$ is a small positive constant. Suppose $\epsilon$ and $\gamma$ are positive constants with $\psi_2(8\epsilon)\le\delta$ and $\gamma\le\delta/(2\log 2)$. Suppose $C'/\gamma \le L \le N\delta/C'$, and divide $[N]$ into consecutive blocks $I_1,\ldots,I_L$ of size $K=N/L$ each. Then there exists $J_\star\subseteq[L]$ with $|J_\star|=N\gamma$,
and $\Omega\subseteq S$ with $|\Omega|\ge \exp( N\delta/L)$, such that
all pairs $\sigma\ne\tau$ in $\Omega$ satisfy
	\[\f{|(\sigma_{I_j},\tau_{I_j})|}{K}
		\le 1- \epsilon\]
for all $j\in J_\star$. That is to say, all pairs $\sigma\ne\tau$ in $\Omega$  are $\epsilon$-separated on $I_j$ for all $j\in J_\star$.
\end{ppn}

The \hyperlink{proof:p.sep.on.blocks.subsample}{proof of Proposition~\ref{p:sep.on.blocks.subsample}} appears at the end of this subsection.

\begin{lem}\label{l:sep.on.blocks}
As in Proposition~\ref{p:sep.on.blocks.subsample}, assume that
$S$ is any subset of $\set{-1,+1}^N$ with $|S|\ge\exp(8N\delta)$;
and that we have $\psi_2(8\epsilon)\le\delta$, $\gamma \le \delta/(2\log2)$, and $C'/\gamma \le L \le N\delta/C'$. Then
	\[
	\mu^{\otimes 2}\bigg(
	\sum_{j\le L} \I\bigg\{
	\f{|(\sigma_{I_j},\tau_{I_j})|}{K}
	\le 1-8\epsilon\bigg\} < 2L\gamma
	\bigg)
	\le \f1{\exp(5N\delta)}
	\,,
	\]
where $\mu$ is the uniform measure on $S$, and $(\sigma,\tau)$ is a pair of i.i.d.\ samples from $\mu$. 

\begin{proof} Write $t\equiv 1-8\epsilon$. For $\sigma\in\set{-1,+1}^N$, we can bound
	\begin{align*}
	\mu_\sigma
	&\equiv
	\mu\bigg(\bigg\{\tau:
	\sum_{j\le L} \I\bigg\{
	\f{|(\sigma_{I_j},\tau_{I_j})|}{K}
	\le 1-8\epsilon\bigg\} < 2L\gamma
	\bigg\}
	\bigg) \\
	&\le \f{1}{|S|}
	\bigg|\bigg\{\tau \in\set{-1,+1}^N:
	\sum_{j\le L} \I\bigg\{
	\f{|(\sigma_{I_j},\tau_{I_j})|}{K}
	\le 1-8\epsilon\bigg\} < 2L\gamma
	\bigg\}\bigg|\,.
	\end{align*}
For any fixed $\sigma$, if $\tau$ is sampled uniformly at random from $\set{-1,+1}^N$, then the scalar product $(\sigma_{I_j},\tau_{I_j})$ is equidistributed as $2(\mathrm{Bin}(K,1/2)-K/2)$, so the chance for 
$|(\sigma_{I_j},\tau_{I_j})| \ge K(1-8\epsilon)$ is 
	\[
	\bar{p}\equiv
	\f1{2^K}\bigg|\bigg\{
	\tau_{I_j}\in\set{-1,+1}^K
	: \f{|(\sigma_{I_j},\tau_{I_j})|}{K} \le t\bigg\}\bigg|
	=\P\bigg( \bigg|\mathrm{Bin}\bigg(K,\f12\bigg)-K
		\bigg|\ge \f{Kt}{2}\bigg)
	\le \f{2}{\exp(K k_2(t))}
	\]
Moreover the $(\sigma_{I_j},\tau_{I_j})$ are independent across the different $j$, so 
	\[\mu_\sigma
	\le \f{2^N}{\exp(8N\delta)}
	\P\bigg( \mathrm{Bin}(L,\bar{p}) \ge L(1-2\gamma)\bigg)  
	\le \f{2^N}{\exp(8N\delta)} \binom{L}{2L\gamma}
	\bar{p}^{L(1-2\gamma)}\,.\]
Combining with the preceding bound on $\bar{p}$, and recalling
$k_2(t)=\log2-\psi_2(8\epsilon)$, we obtain
	\begin{align*}
	\mu_\sigma
	&\le\f{2^N \exp(LH(2\gamma))}{\exp(8N\delta)}
	\bigg(
	\f{2}{\exp(K k_2(t))}
	\bigg)^{L(1-2\gamma)} \\
	&\le \exp\bigg\{
	N\Big( -8\delta + \psi_2(8\epsilon)
	+2\gamma\log2 - 2\gamma\psi_2(8\epsilon)
	\Big)
	+ L\Big( H(2\gamma) + \log 2 \Big)
	\bigg\}\le \f1{\exp(5N\delta)}\,,
	\end{align*}
where the last bound holds for parameters $\delta,\epsilon,\gamma,L$ as in the statement of the lemma. The result follows by averaging over $\sigma\in S$.
\end{proof}
\end{lem}

\begin{cor}\label{c:sep.on.blocks}
In the setting of Lemma~\ref{l:sep.on.blocks} (and with all the same parameters),
there must be a subset
$S'\subseteq S$
with $|S'|\ge \exp(5N\delta/2)$
 such that all pairs $\sigma\ne\tau$ in $S'$ 
 are well-separated in the sense that
	\beq\label{e:well.sep}
	\sum_{j\le L} \I\bigg\{
	\f{|(\sigma_{I_j},\tau_{I_j})|}{K}
	\le 1-8\epsilon\bigg\} \ge 2L\gamma\,.\eeq

\begin{proof}
Let $\sigma^{(1)},\ldots,\sigma^{(n)}$ be $n=\exp(5N\delta/2)$ i.i.d.\ samples from $\mu$ (the uniform measure on $S$). By Lemma~\ref{l:sep.on.blocks}, the probability that for any $1\le i<j\le n$ the configurations $\sigma^i,\sigma^j$ fail to be well-separated is upper bounded by
	\[
	\binom{n}{2}
	\mu^{\otimes 2}
	\bigg(
	\sum_{j\le L} \I\bigg\{
	\f{|(\sigma_{I_j},\tau_{I_j})|}{K}
	\le 1-8\epsilon\bigg\} < 2L\gamma
	\bigg)
	\le
	\binom{n}{2}
	\f1{\exp(5N\delta)}
	\le \f12\,,
	\]
where the last bound holds by the choice of $n$. Therefore the random set $\set{\sigma^{(1)},\ldots,\sigma^{(n)}}$ satisfies the required condition with probability at least $1/2$; and this implies the existence of the claimed set $S'$.
\end{proof}
\end{cor}

\begin{lem}\label{l:sep.on.fixed.blocks}Let $S'\subseteq \set{-1,+1}^N$ with $|S'|\ge \exp(5N\delta/2)$, such that all pairs $\sigma\ne\tau$ in $S'$ are well-separated in the sense of  \eqref{e:well.sep}. If $\eta$ satisfies
$(2/|S'|)^{1/L} \le \eta\le1/2$, then there exist $J_\circ\subseteq [L]$ with $|J_\circ|\ge 2L\gamma$, and $T\subseteq S'$ with $|T| \ge \eta^L|S'|$, such that we have
	\beq\label{e:sep.bound}
	\f1{|T|}\bigg|\bigg\{\tau\in T
		: \f{|(\sigma_{I_j},\tau_{I_j})|}{K}
			> 1- \epsilon
		\bigg\}
		\bigg| \le 2\eta\,.
	\eeq
for all $\sigma\in T$ and all $j\in J_\circ$.

\begin{proof}
We first give the construction of $T$ and $J_\circ$. For $\sigma\in\set{-1,+1}^N$, $\Omega\subseteq\set{-1,+1}^N$, and $1\le j\le L$, denote
	\[T_j(\sigma,\Omega)
	\equiv\bigg\{\tau\in\Omega:
		\f{(\sigma_{I_j},\tau_{I_j})}{K} > 1-\epsilon\bigg\}\,.\]
\begin{enumerate}[$\bullet$]
\item \textit{Initial step.}
Let $\Omega^{(0)}\equiv S'$. If for all $\sigma\in\Omega^{(0)}$ and all $j\in[L]$ we have
	\[
	\max\bigg\{
	\f{|T_j(\sigma,\Omega^{(0)})|}{|\Omega^{(0)}|},
	\f{|T_j(-\sigma,\Omega^{(0)})|}{|\Omega^{(0)}|}\bigg\}\le \eta\,,
	\]
then we are done by simply taking $T = \Omega^{(0)}$ and $J_\circ=[L]$. If this is not the case, then we must have a configuration $\sigma^{(1)}\in\Omega^{(0)} \cup (-\Omega^{(0)})$ and an index $j^{(1)}\in[L]$ such that
	\[\f{|T_{j^{(1)}}(\sigma^{(1)},\Omega^{(0)})|}{|\Omega^{(0)}|}
	> \eta\,.
	\]
Then set $\Omega^{(1)}=T_{j^{(1)}}(\sigma^{(1)},\Omega^{(0)})\subseteq\Omega^{(0)}$.
\item \textit{Inductive step.} Similarly, suppose inductively that we have constructed
$\sigma^{(k)}$, $j^{(k)}$, $\Omega^{(k)}$. If for all $\sigma\in\Omega^{(k)}$ 
and all $j\in[L]\setminus\set{j^{(1)},\ldots,j^{(k)}}$
we have
	\[
	\max\bigg\{
	\f{|T_j(\sigma,\Omega^{(k)})|}{|\Omega^{(k)}|},
	\f{|T_j(-\sigma,\Omega^{(k)})|}{|\Omega^{(k)}|}\bigg\}\le \eta\,,
	\]
then we end the process
by setting $T=\Omega^{(k)}$  and $J_\circ=[L]\setminus\set{j^{(1)},\ldots,j^{(k)}}$. If this is not the case, then we must have a configuration $\sigma^{(k+1)}\in\Omega^{(k)} \cup (-\Omega^{(k)})$ and an index $j^{(k+1)}\in[L]\setminus\set{j^{(1)},\ldots,j^{(k)}}$ such that
	\[\f{|T_{j^{(k+1)}}(\sigma^{(k+1)},\Omega^{(k)})|}{|\Omega^{(k)}|}
	> \eta\,.
	\]
Then set $\Omega^{(k+1)}
=T_{j^{(k+1)}}(\sigma^{(k+1)},\Omega^{(k)})\subseteq\Omega^{(k)}$.
\end{enumerate}
We claim that the above construction ends with $T=\Omega^{(k')}$ for some $k' \le \ell = L(1-2\gamma)$. Indeed, suppose for contradiction that it does not. For simplicity of notation, let us re-index such that $j^{(k)}=k$ for all $k$. Now consider the set $F=
\Omega^{(\ell+1)}=
T_{\ell+1}(\sigma^{(\ell+1)},\Omega^{(\ell)})$
for $\ell=L(1-2\gamma)$. If $\tau\in F$, then for each $k\le\ell+1$ the configurations $\tau$ and $\sigma^{(k)}$ must be close on the block of coordinates $I_k$:
	\[
	\f{\|\tau_{I_k}-(\sigma^{(k)})_{I_k}\|^2}{K}
	= \f{2K-2(\tau_{I_k},(\sigma^{(k)})_{I_k})}{K}
	\le 2-2(1-\epsilon) = 2\epsilon\,.
	\]
By the assumption $\eta^L|S'| \ge 2$ we must have
	\[
	|F| 
	= |\Omega^{(0)}| \prod_{k=1}^{\ell+1}
	\f{|\Omega^{(k)}|}{|\Omega^{(k-1)}|}
	\ge |S'| \eta^{\ell+1}
	\ge |S'| \eta^L \ge 2\,,
	\]
so we can find a pair of distinct elements $\sigma\ne\tau$ in $F$. For this pair, the triangle inequality gives
	\[
	\f{\|\sigma_{I_k}-\tau_{I_k}\|^2}{K}
	\le \Big(
	2 (2\epsilon)^{1/2}\Big)^2
	 =8 \epsilon\]
for each $k\le \ell+1$. This contradicts the assumption that $\sigma$ and $\tau$ must satisfy \eqref{e:well.sep}, since they are both elements of the original set $S'$. This verifies the claim that $T=\Omega^{(k')}$ for some $k'\le \ell = L(1-2\gamma)$, and 
	\[
	\f{|T|}{|S'|}
	= \prod_{k=1}^{k'}
	\f{|\Omega^{(k)}|}{|\Omega^{(k-1)}|}
	\ge \eta^{k'} \ge \eta^\ell
	\ge \eta^L\,.
	\]
Finally, the required bound \eqref{e:sep.bound} holds because otherwise the above construction would not stop at $T=\Omega^{(k')}$.
\end{proof}
\end{lem}

\begin{lem}\label{l:all.pairs.sep.on.fixed.blocks}
Let $S$ be any subset of $\set{-1,+1}^N$ with $|S|\ge\exp(8N\delta)$ where $\delta$ is a small positive constant. Suppose $\psi_2(8\epsilon)\le\delta$, $\gamma\le\delta/(2\log 2)$, 
and $C'/\gamma \le L \le N\delta/C'$. Then there exist a subset $J_\circ\subseteq[L]$ with $|J_\circ|\ge 2L\gamma$, and a subset $T\subseteq S$ with $|T| \ge \exp(N\delta)$, such that for all $\sigma\in T$ we have
	\[
	\f1{|T|}\bigg|\bigg\{\tau\in T
		: \f{|(\sigma_{I_j},\tau_{I_j})|}{K}
			> 1- \epsilon
		\bigg\}
		\bigg| \le 2\exp\bigg(-\f{9N\delta}{4L}\bigg)
	\]
for each $j\in J_\circ$. 

\begin{proof}
Let $S$ be any subset of $\set{-1,+1}^N$ with $|S|\ge \exp(8N\delta)$. It follows from Corollary~\ref{c:sep.on.blocks} that there must be a subset $S'\subseteq S$ with $|S'|\ge \exp(5N\delta/2)$ such that all pairs $\sigma\ne\tau$ in $S'$ are well-separated in the sense of \eqref{e:well.sep}. It then  follows from Lemma~\ref{l:sep.on.fixed.blocks} that if $(2/|S'|)^{1/L} \le \eta\le1/2$, then there exists $J_\circ\subseteq[L]$ with $|J_\circ|\ge 2L\gamma$,
and  $T\subseteq S'$
with $|T| \ge \eta^L|S'|$, such that we have
	\[
	\f1{|T|}\bigg|\bigg\{\tau\in T
		: \f{|(\sigma_{I_j},\tau_{I_j})|}{K}
			> 1- \epsilon
		\bigg\}
		\bigg| \le 2\eta
	\]
for all $\sigma\in T$ and all $j\in J_\circ$. In particular, we can choose $\eta=\exp(-9N\delta/(4L))$, so that 
	\[\eta^L
	= \exp\bigg(-\f{9N\delta}{4}\bigg) 
	\ge 2\exp\bigg(-\f{5N\delta}{2}\bigg)
	\ge \f{2}{|S'|}
	\]
holds for all $N$ large enough, as required.
\end{proof}
\end{lem}

\begin{proof}[\hypertarget{proof:p.sep.on.blocks.subsample}{Proof of Proposition~\ref{p:sep.on.blocks.subsample}}] 
Let $T\subseteq S$ be as given by Lemma~\ref{l:all.pairs.sep.on.fixed.blocks}, and let $\mu$ be the uniform probability measure on $T$. Similarly to the proof of Corollary~\ref{c:sep.on.blocks}, let $(\sigma^{(1)},\ldots,\sigma^{(n)})$ be sampled according to $\mu^{\otimes n}$, i.e., uniformly at random from $T^n$. Let us say that a block $I_k$ ``fails'' if there exists a pair of indices $1\le i<j\le n$ such that the configurations $\sigma^{(i)}$ and $\sigma^{(j)}$ are not $\epsilon$-separated on the block $I_k$. 
 The probability that more than half of the blocks in $J_\circ$ fail can be upper bounded with Markov's inequality:
	\begin{align*}
	\mu^{\otimes n}\bigg(
	\f1{|J_\circ|}
	\sum_{j\in J_\circ}
	\Ind{\textup{block $I_j$ fails}}\ge\f12
	\bigg)
	&\le \f{2}{|J_\circ|}
	\sum_{j\in J_\circ}
	\binom{n}{2}
	\mu^{\otimes 2}\bigg(\f{|(\sigma_{I_j},\tau_{I_j})|}{K}
		> 1- \epsilon\bigg) \\
	&\le n^2 \exp\bigg( -\f{9N\delta}{4L} \bigg) \le\f12\,,
	\end{align*}
where the last inequality holds for $n\le \exp(N\delta/L)$.
Therefore, with probability at least $1/2$,
for the random set 
$\set{\sigma^{(1)},\ldots,\sigma^{(n)}}$
we can find $|J_\star|\ge L\gamma$ such that the required properties are satisfied. If we take $n= \exp(N\delta/L)$, then this implies the existence  of the claimed set $\Omega\subseteq T\subseteq S$.
\end{proof}

\subsection{Intersection of cube and half-space with general disorder}
\label{ss:threshold.general}

The main results in this subsection are Proposition~\ref{p:largevals.general},
which is a version of
Proposition~\ref{p:largevals.gaus} for general disorder,
and the
\hyperlink{proof:t.hspace.general}{proof of Theorem~\ref{t:hspace.general}}.
The assumption of the following is based on the result of Proposition~\ref{p:sep.on.blocks.subsample}:

\begin{ppn}\label{p:general.disorder.threshold}
Suppose $\XI$ is a random vector in $\R^N$ satisfying Assumption~\ref{a:subgaus}. Suppose $\delta$ is a small positive constant, and $\epsilon,\gamma\le\delta$. Divide $[N]$ into consecutive blocks $I_1,\ldots,I_L$ of size $K=N/L$ each. Suppose $\Omega\subseteq\set{-1,+1}^N$
with $|\Omega|\ge \exp(N\delta/L)$, such that
all pairs $\sigma\ne\tau$ are $\epsilon$-separated on $I_j$ for all 
for all $j\in J_\star = [\ell+1,L]$, where $\ell=L(1-\gamma)$.  Then we have the bound 
	\[
	\P
	\bigg(
	\f1{|\Omega|}
	\bigg|\bigg\{\sigma\in\Omega:
		\f{(\XI,\sigma)}{N^{1/2}}
		\ge \f{\gamma(L\epsilon\log p)^{1/2}}{60}
		\bigg\}\bigg| 
		< \f1{2\cdot(4p)^{L\gamma}}
	\bigg)
	\le 
	4\exp\bigg(-\f{L\gamma^2\epsilon\log p}{60\NU}\bigg)\,,
	\]
provided $\log p \ge C'\NU/\epsilon$
and $C'/\gamma \le L \le N^{1/2}/(C'\log p)^{1/2}$.

\begin{proof}
Let $p$ be a large integer. 
Let $\filt_k$ be the $\sigma$-algebra generated by the vectors $\XI_{I_j}$ for $1\le j\le k$, so that $M_k(\sigma)$
(as defined by \eqref{e:block.process}) is measurable with respect to $\filt_k$. 
Recall $\ell=L(1-\gamma)$, and define
	\[
	\Omega^{(\ell)}
	\equiv\bigg\{ \sigma\in\Omega :
		M_\ell(\sigma) \ge
		- \f{\gamma(L\epsilon\log p)^{1/2}}{5}
		\bigg\}\,.
	\]
It follows using Markov's inequality and the subgaussian tail bound that
	\begin{align*}
	\textup{(I)} &\equiv
	\P
	\bigg(
	\f{|\Omega^{(\ell)}|}{|\Omega|}
	\le \f12\bigg)
	\le \P
	\bigg(
	\f1{|\Omega|} \bigg|\bigg\{\sigma\in\Omega: 
		M_\ell(\sigma) <
		- \f{\gamma(L\epsilon\log p)^{1/2}}{5}
		\bigg\}\bigg|
	\ge \f12
	\bigg)\\
	&\le 
	2\exp\bigg( -\f{L \gamma^2 \epsilon\log p}{50\NU} \bigg) \,.
	\end{align*}
Next, for $\ell+1\le k\le L$, suppose inductively that $\Omega^{(k)}$ has been defined. Let
	\[
	\Xi^{(k)}
	\equiv 
	\bigg\{
	\sigma\in\Omega^{(k)}
	: \f{(\XI_{I_k},\sigma_{I_k})}{K^{1/2}}
	\ge \f{(\epsilon\log p)^{1/2}}{2}
	\bigg\}\subseteq\Omega^{(k)}\,.
	\]
We will show below that $\Xi^{(k)}$ occupies a constant fraction of $\Omega^{(k)}$ with good probability:
	\beq\label{e:Omega.Xi.general.disorder}
	\P\bigg(
	\f{|\Xi^{(k)}|}{|\Omega^{(k)}|}
	\le \f1{4p}
	\bigg) \le \f{4}{p^{\epsilon/32}}\,.
	\eeq
With this in mind, for $\ell+1\le k\le L$ we define the sets
	\[
	\Omega^{(k+1)}=\begin{cases}
	\Xi^{(k)} & \textup{if
		$|\Xi^{(k)}|
		\ge |\Omega^{(k)}| / (4p)$}\\
	\Omega^{(k)} & \textup{otherwise.}
	\end{cases}
	\]
The construction guarantees for all $\ell \le k\le L$ that
	\beq\label{e:lbd.p}
	|\Omega^{(k)}|
	\ge |\Omega^{(L)}|
	\ge \f{|\Omega|}2 \bigg(\f1{4p}\bigg)^{L\gamma}
	\ge \f12 \exp\bigg( \f{N\delta}{L} - L\gamma\log(4p)\bigg)
	\ge p\,,\eeq
where the last bound holds by the assumed upper bound on $L$. Now, turning to the proof of \eqref{e:Omega.Xi.general.disorder}, from $\Omega^{(k)}$ let us extract disjoint subsets $X_1,\ldots,X_m$ each of size $p$, whose union occupies more than half of $\Omega^{(k)}$ --- the preceding bound guarantees that we can do this with $m\ge1$, since $|\Omega^{(k)}|\ge p$. Combining Corollary~\ref{c:all.fail} with the CLT estimate Corollary~\ref{c:clt} (deferred to \S\ref{ss:clt} below) gives
	\begin{align*}
	\P\bigg( \max_{\sigma\in X_a}
		\f{(\XI,\sigma)}{N^{1/2}}
		\le \f{(\epsilon\log p)^{1/2}}{2}\bigg)
	&= \E\bigg[ \prod_{\sigma\in X_a}
	\I\bigg\{
	\f{(\XI,\sigma)}{N^{1/2}}  \le \f{(\epsilon\log p)^{1/2}}{2}
	\bigg\}
	\bigg] \\
	&\le
	\P\bigg( \max_{\sigma\in X_a} \f{(g,\sigma)}{N^{1/2}}
		\le \f{(\epsilon\log p)^{1/2}}{2}\bigg) + o_N(1)
	\le \f{2}{p^{\epsilon/50}}
	\end{align*}
It then follows by the Markov inequality that
	\[
	\P\bigg(
	\f1m\sum_{a\le m} \I\Big\{
	\Xi^{(k)} \cap X_a=\varnothing\Big\}
	\ge \f12
	\bigg)
	\le \f{4}{p^{\epsilon/50}}\,.
	\]
On the complementary event we have
	\[
	|\Xi^{(k)}|
	\ge \sum_{a\le m}
	\Ind{\Xi^{(k)} \cap X_a\ne\varnothing}
	\ge \f{m}{2}
	\ge \f{|\Omega^{(k)}|}{4p}\,,
	\]
which proves \eqref{e:Omega.Xi.general.disorder}.
Next let $J_\bullet$ denote the subset of indices $k\in J_\star$ for which we have $\Omega^{(k)}=\Omega^{(k-1)}$. It follows from \eqref{e:Omega.Xi.general.disorder} that $|J_\bullet|$ is stochastically dominated by a binomial random variable with $L\gamma$ trials and success probability $4/p^{\epsilon/50}$. It follows by the Chernoff bound and \eqref{e:chernoff.simplified} that
	\[
	\textup{(II)}
	\equiv
	\P\bigg(|J_\bullet| 
	\ge \f{L\gamma}{6\NU}
	\bigg)
	\le
	\exp\bigg\{ -L\gamma H\bigg( 
		\f1{6\NU}
		\bigg| \f{4}{p^{\epsilon/50}} \bigg)\bigg\}
	\le
	\exp\bigg( -\f{L\gamma \epsilon\log p}
	{6\NU \cdot 51} \bigg)
	\,,
	\]
where the last bound holds because the restriction $\log p\ge C'\NU/\epsilon$ guarantees
	\[
	\f{t}{e} = \f1{e} \cdot \f1{6\NU}
		\cdot \f{p^{\epsilon/50}}{4} \ge p^{\epsilon/51}\,.
	\]
By a union bound over all subsets $J\subseteq J_\star$, we have
	\begin{align*}
	\textup{(III)}
	&\equiv
	\P\bigg(
	\sum_{j\in J} \f{(\XI_{I_j},\sigma_{I_j})}{N^{1/2}}
	\le -\f{\gamma(L\epsilon\log p)^{1/2}}{5}
	\textup{ for any }J\subseteq J_\star
	\bigg) \\
	&\le
	2^{L\gamma}
	\exp\bigg( -\f{L\gamma\epsilon\log p}{50\NU} \bigg)
	\le 
	\exp\bigg( -\f{L\gamma\epsilon \log p}{60\NU} \bigg)\,,
	\end{align*}
where the last bound again uses the restriction $\log p\ge C'\NU/\epsilon$. On the complement of the events bounded by $\textup{(I)}$, $\textup{(II)}$, and $\textup{(III)}$, we have for all $\sigma\in\Omega^{(L)}$ that
	\begin{align*}
	\f{(\XI,\sigma)}{N^{1/2}}
	&= M_L(\sigma)
	= M_\ell(\sigma)
	+ \sum_{j\in J} 
	\f{(\XI_{I_j},\sigma_{I_j})}{N^{1/2}}
	+ \sum_{j\in J_\star\setminus J_\bullet} 
	\f{(\XI_{I_j},\sigma_{I_j})}{N^{1/2}} \\
	&\ge
	- \f{2\gamma(L\epsilon\log p)^{1/2}}{5}
	+ L\gamma\bigg(1-\f1{6\NU}\bigg)
	\f{(\epsilon\log p)^{1/2}}{2 L^{1/2}}
	\ge \f{\gamma(L\epsilon\log p)^{1/2}}{60}\,,
	\end{align*}
where the last bound uses that we must have $\NU\ge1$ (see Assumption~\ref{a:subgaus}). The claim then follows by recalling the lower bound on $|\Omega^{(L)}|$ from \eqref{e:lbd.p}.
\end{proof}
\end{ppn}

The next proposition is a version of Proposition~\ref{p:largevals.gaus}
which applies in the case of general disorder. It has a somewhat worse $\epsilon$-dependence Proposition~\ref{p:largevals.gaus}. From our perspective the more important difference is that 
Proposition~\ref{p:largevals.general} applies for a more limited range of $s$ than Proposition~\ref{p:largevals.gaus}:

\begin{ppn}\label{p:largevals.general}
Let $S$ be any subset of $\set{-1,+1}^N$ with $|S|\ge\exp(N\delta)$ where $\delta$ is a small positive constant. Let $\epsilon$ be a positive constant with $\psi_2(8\epsilon)\le \delta/9$. Then 
	\[\P\bigg(
	\f1{|S|}
		\bigg|\bigg\{\sigma\in S: 
	\f{(\XI,\sigma)}{N^{1/2}}
		\ge s(\epsilon\NU)^{1/2} 
		\bigg\}\bigg|
	< \exp\bigg( -\f{C''\NU s^2}{2\epsilon}\bigg)
	\bigg)
	\le\f{4}{\exp(s^2\epsilon)}\,.
	\]
provided  $C'' \le s \le (N\epsilon^3/\NU)^{1/4}/C''$, where $C''$ is a large absolute constant.

\begin{proof}
Since we can relabel the constants later, we assume $|S|\ge\exp(9N\delta)$ and $\psi_2(8\epsilon)\le\delta$. Recalling the statement of Proposition~\ref{p:sep.on.blocks.subsample}, let $\gamma=\epsilon$. (For $\delta$ small enough, the assumption $\psi_2(8\epsilon)\le\delta$ and the choice $\gamma=\epsilon$ together guarantee that the condition $\gamma\le \delta/(2\log2)$ from Proposition~\ref{p:sep.on.blocks.subsample} is satisfied.) It follows from Proposition~\ref{p:sep.on.blocks.subsample} that there exists $J_1\subseteq[L]$ with $|J_1|=L\gamma$, and $\Omega_1\subseteq S$ with $|\Omega_1|=\exp(N\delta/L)$, such that all pairs $\sigma\ne\tau$ in $\Omega_1$ are $\epsilon$-separated on $I_j$ for all $j\in J_1$. We can apply Proposition~\ref{p:sep.on.blocks.subsample} again on $S\setminus\Omega_1$, and so on, to extract
disjoint subsets $\Omega_1,\ldots,\Omega_m\subseteq S$, up to the first $m$ such that
	\[
	\sum_{a\le m}|\Omega_a| \ge \f{|S|}{2}\,.
	\]
Then for each $a$ we will have $|\Omega_a|=\exp(N\delta/L)$,
and all pairs $\sigma\ne\tau$ in $\Omega_a$ will be $\epsilon$-separated on $I_j$ for all $j\in J_a$, where $J_a\subseteq[L]$ with $|J_a|=L\gamma$. 
Recalling the statement of Proposition~\ref{p:general.disorder.threshold}, 
let $\log p=C'\NU/\epsilon$, and define
	\[
	\Xi\equiv\bigg\{\sigma\in S: 
		\f{(\XI,\sigma)}{N^{1/2}}
		\ge \f{\gamma(L\epsilon\log p)^{1/2}}{60} 
		\bigg\}\,,\quad
	\Xi_a
	\equiv \Xi \cap \Omega_a\,.
	\]
By Markov's inequality and Proposition~\ref{p:general.disorder.threshold},
as long as $C'/\gamma \le L \le N^{1/2}/(C'\log p)^{1/2}$, we have 
	\[
	\P\bigg(\f1m \sum_{a\le m}
	\I\bigg\{ \f{|\Xi_a|}{|\Omega_a|}
	\le \f1{2\cdot (4p)^{L\gamma}}
	\bigg\} \ge \f12\bigg)
	\le 4\exp\bigg(-\f{L\gamma^2\epsilon\log p}{60\NU}\bigg)\,.
	\]
On the complementary event we must have
	\[
	|\Xi|
	\ge \sum_{a\le m} |\Xi_a|
	\ge \f{m|\Omega_1|}{2} 
	\cdot \f1{2\cdot (4p)^{L\gamma}}
	\ge \f{|S|}{8\cdot(4p)^{L\gamma}}
	\,,
	\]
so we have shown that, for any $L$ satisfying 
$C'/\gamma \le L \le N^{1/2}/(C'\log p)^{1/2}$, we have
	\[
	\P\bigg(
	\f1{|S|}
		\bigg|\bigg\{\sigma\in S: 
	\f{(\XI,\sigma)}{N^{1/2}}
		\ge \f{\gamma(L\epsilon\log p)^{1/2}}{60}
		\bigg\}\bigg|
	< \f1{8\cdot(4p)^{L\gamma}}
	\bigg)
	\le
	4\exp\bigg(-\f{L\gamma^2\epsilon\log p}{60\NU}\bigg)\,.
	\]
Recalling that $\log p=C'\NU/\epsilon$, the above implies
	\[
	\P\bigg(
	\f1{|S|}
		\bigg|\bigg\{\sigma\in S: 
	\f{(\XI,\sigma)}{N^{1/2}}
		\ge \f{(C'L\NU)^{1/2}\gamma}{60}
		\bigg\}\bigg|
	<\f1{8} \exp\bigg( -\f{2 C'L\NU \gamma}{\epsilon}
	\bigg)
	\bigg)
	\le
	4\exp\bigg(-\f{C'L\gamma^2}{60}\bigg)\,.
	\]
Making the change of variables $s=(C'L\gamma)^{1/2}/60$ gives
	\[
	\P\bigg(
	\f1{|S|}
		\bigg|\bigg\{\sigma\in S: 
	\f{(\XI,\sigma)}{N^{1/2}}
		\ge s(\gamma\NU)^{1/2}
		\bigg\}\bigg|
	< \f1{8}
	\exp\bigg( -\f{(2\cdot 60^2) \NU s^2}{\epsilon}
	\bigg)
	\bigg)
	\le
	\f{4}{\exp(60 \gamma s^2)}\,.
	\]
and the conclusion follows by recalling that we took $\gamma=\epsilon$. Note that the bounds on $s$ in the statement of the result guarantee that $L$ satisfies the requirements 
$C'/\gamma \le L \le N^{1/2}/(C'\log p)^{1/2}$.
\end{proof}
\end{ppn}

\begin{proof}[\hypertarget{proof:t.hspace.general}{Proof of Theorem~\ref{t:hspace.general}}]
Let $\epsilon$ be a positive constant with $\psi_2(8\epsilon)\le \delta/9$. It follows directly from Proposition~\ref{p:largevals.general} that
	\[\P\bigg(
	\f1{|S|}
		\bigg|\bigg\{\sigma\in S: 
	\f{(\XI,\sigma)}{N^{1/2}}
		\ge \kappa\bigg\}\bigg|
	< \exp\bigg( -\f{C''\NU s^2}{2\epsilon}\bigg)
	\bigg)
	\le\f{4}{\exp(s^2\epsilon)}\,,
	\]
provided that $s$ satisfies (with $C''$ a large absolute constant) 
	\[
	\max\bigg\{
	C'', \f{\kappa_+}{(\epsilon\NU)^{1/2}}
	\bigg\} \le s
	\le\f1{C''} \bigg( \f{N\epsilon^3}{\NU}\bigg)^{1/4}\,.
	\]
The claim follows by fixing $\epsilon$ to satisfy $\psi_2(8\epsilon)= \delta/9$, and making the change of variables $w= C''\NU s^2/(2\epsilon)$.
\end{proof}

\subsection{Intersection of cube and slab with general disorder} 
\label{ss:addone.general}

The \hyperlink{proof:t.interval.general}{proof of Theorem~\ref{t:interval.general}} appears at the end of this subsection. 
The next two lemmas are motivated by the following considerations. In view of Proposition~\ref{p:sep.on.blocks.subsample}, suppose we have a large subset $\Omega\subseteq\set{-1,+1}^N$ such that all pairs $\sigma\ne\tau$ in $\Omega$ are well-separated on the blocks $I_j$ for $j\in J_\star$. By re-indexing, we may suppose without loss that $J_\star = [L]\setminus [\ell]$ for $\ell=L(1-\gamma)$. For $\sigma\in\Omega$, recall the process defined by \eqref{e:block.process}, and suppose we have $|M_\ell(\sigma)|\le s$ for $s$ large. Then, in order to achieve the desired outcome $M_L(\sigma)\in[a,b]$, the process $M_k(\sigma)$ may need to traverse a distance at most $2s$ over the final $L\gamma$ blocks, $k\in[\ell+1,L]$. This suggests that we consider the events
	\[
	\bigg\{\f{(g_{I_k},\sigma_{I_k})}{N^{1/2}} \ge \f{2s}{L\gamma}\bigg\}
	\,,\quad
	\bigg\{-\f{(g_{I_k},\sigma_{I_k})}{N^{1/2}} \ge \f{2s}{L\gamma}\bigg\}
	\]
so that the process can traverse enough distance on each block. On the other hand, if $M_{k-1}(\sigma)\in[a,b]$, then we can ensure $M_k(\sigma)\in[a,b]$ by requiring one of the two events
	\[
	\bigg\{0\le\f{(g_{I_k},\sigma_{I_k})}{N^{1/2}} \le \f{b-a}{2}\bigg\}\,,\quad
	\bigg\{0\le -\f{(g_{I_k},\sigma_{I_k})}{N^{1/2}} \le \f{b-a}{2}\bigg\}\,.
	\]
This leads to the statement of the following:

\begin{lem}\label{l:Kbound.gaus} Let $\chi\in(0,1]$ be a constant. Let $S$ be any subset of $\set{-1,+1}^K$ with $|S|\ge\exp(K\delta)$ where $\delta$ is a positive constant. Let $\epsilon$ be a positive constant with $\psi_2(\epsilon)\le\delta/2$.
	\[
	\P\bigg(
	\f1{|S|}
	\bigg|
	\bigg\{\sigma\in S : \f{(g,\sigma)}{N^{1/2}} 
	\in\bigg[ \f{r\epsilon^{1/2}}{L^{1/2}} ,\chi
	\bigg]
	\bigg\}\bigg|
	\le \f1{\exp(2C' s^2)}
	\bigg) 
	\le \f2{\exp(s^2\epsilon/C')}\,,
	\]
provided that the parameters satisfy the bounds 
	\[
	C' \le r\le s \le 
	\min\bigg\{
	\f{L^{1/2}\chi}{(6C')^{1/2}}
	, \f{(K\delta)^{1/2}}{C'}
	\bigg\}\,.
	\]
\begin{proof}
It follows directly from Proposition~\ref{p:largevals.gaus} that for $C' \le r\le s \le (K\delta)^{1/2} /C'$ we have
	\[
	\P\bigg(
	\f1{|S|}
	\bigg|
	\bigg\{\sigma\in S : \f{(g,\sigma)}{N^{1/2}} \ge 
	\f{r\epsilon^{1/2}}{L^{1/2}} 
	\bigg\}\bigg|
	\le \f1{\exp(C' s^2)}
	\bigg) 
	\le \f{1}{\exp(s^2\epsilon/C')}\,,
	\]
On the other hand, it following using Markov's inequality and the gaussian tail bound that
	\[\P\bigg(\f1{|S|}
	\bigg|\bigg\{ \sigma\in S:
	\f{(g,\sigma)}{N^{1/2}} \ge \chi \bigg\}\bigg|
	\ge  \f1{\exp(2C's^2)}
	\bigg)
	\le
	\f{\exp(2 C' s^2)}
		{\exp(\chi^2 L/2)}
	\le 
	\f1{\exp(C's^2)}\,,
	\]
where the last bound holds provided that $s^2 \le L\chi^2/(6C')$. Combining these bounds gives the claim.
\end{proof}
\end{lem}

The following is a version of Lemma~\ref{l:Kbound.gaus} for general disorder. It has a worse $\epsilon$-dependence and applies in a more limited range of $s$, as a consequence of applying Proposition~\ref{p:largevals.general} in place of Proposition~\ref{p:largevals.gaus}.

\begin{lem}\label{l:Kbound.general} Let $\chi\in(0,1]$ be a constant. Let $S$ be any subset of $\set{-1,+1}^K$ with $|S|\ge\exp(K\delta)$ where $\delta$ is a small positive constant. Let $\epsilon$ be a positive constant with $\psi_2(8\epsilon)\le\delta/9$. Then
	\[\P\bigg(
	\f1{|S|}
		\bigg|\bigg\{\sigma\in S: 
	\f{(\XI,\sigma)}{N^{1/2}}
		\in
		\bigg[\f{r(\epsilon\NU)^{1/2}}{L^{1/2}},
		 \chi\bigg]
		\bigg\}\bigg|
	< \exp\bigg( -\f{C''\NU s^2}{\epsilon}\bigg)
	\bigg)
	\le\f{5}{\exp(s^2\epsilon)}\,,
	\]
provided that the parameters satisfy the bounds 
	\[
	C'' \le r\le s \le 
	\min\bigg\{
	\f{(L\epsilon)^{1/2}\chi}{2(C'')^{1/2}\NU},
	\f1{C''} \bigg( \f{K\epsilon^3}{\NU} \bigg)^{1/4}
	\bigg\}\,.
	\]

\begin{proof}
It follows directly from Proposition~\ref{p:largevals.general} that for $C'' \le r\le s \le (K\epsilon^3/\NU)^{1/4}/C''$ we have
	\[\P\bigg(
	\f1{|S|}
		\bigg|\bigg\{\sigma\in S: 
	\f{(\XI,\sigma)}{N^{1/2}}
		\ge\f{r(\epsilon\NU)^{1/2}}{L^{1/2}}
		\bigg\}\bigg|
	<\exp\bigg( -\f{C''\NU s^2}{2\epsilon}\bigg)
	\bigg)
	\le\f{4}{\exp(s^2\epsilon)}\,.
	\]
On the other hand, it following using Markov's inequality and the subgaussian tail bound that
	\[\P\bigg(\f1n
	\bigg|\bigg\{ i\le n:
	\f{u_i}{L^{1/2}} \ge \chi \bigg\}\bigg|
	\ge  \f1{\exp(C''\NU s^2/\epsilon)}
	\bigg)
	\le
	\f{\exp(C''\NU s^2/\epsilon)}{\exp(\chi^2 L/(2\NU))}
	\le 
	\exp\bigg( -\f{C''\NU s^2}{\epsilon} \bigg)  
	\,,
	\]
where the last bound holds provided that 
$s^2 \le L\chi^2\epsilon/(4C''\NU^2)$.
Combining these bounds gives the claim.
\end{proof}
\end{lem}

The assumption of the next proposition is based on the result of Proposition~\ref{p:sep.on.blocks.subsample}:

\begin{ppn}\label{p:addinterval.general}
Let $\chi\equiv\min\set{1,(b-a)/2}$. Define
	\beq\label{e:defn.RR}
	\RR^2\equiv\f{\rr}{\gamma}
	\equiv
	\f{\max\set{|a|,|b|,C'}}{\gamma}\,.
	\eeq
As before, divide $[N]$ into consecutive blocks $I_1,\ldots,I_L$ of size $K=N/L$ each, such that $L$ satisfies
	\beq\label{e:restrict.L.general.disorder}
	\f{4C''\NU^2\bar{r}}{(\chi\epsilon)^2}
	\le L 
	\le \f{(N\epsilon)^{1/3}\NU}{C''\chi^2}
	\,.\eeq
Now suppose $\Omega\subseteq\set{-1,+1}^N$
with $|\Omega|\ge \exp(2N\delta/L)$, such that
all pairs $\sigma\ne\tau$ are $\epsilon$-separated on $I_j$ for all 
for all $j\in J_\star = [\ell+1,L]$, where $\ell=L(1-\gamma)$.  Then we have the bound 
	\[\P\bigg( \f1{|\Omega|}
	\bigg|\bigg\{\sigma\in\Omega:
	\f{(\XI,\sigma)}{N^{1/2}}
		\in[a,b]\bigg\}\bigg|
	\le
	\f14
	\exp\bigg(-\f{8C''\NU \bar{r} L}{\epsilon}\bigg)
	 \bigg)\le
	\exp\bigg( -\f{L\epsilon^2}{2} \bigg)\,,
	\]
provided that $\gamma=\epsilon\le\delta$, and $\delta$ is small enough.

\begin{proof} 
We emphasize that all the randomness is in the $N$-dimensional random vector $\XI$. Let $\filt_k$ be the $\sigma$-algebra generated by the vectors $\XI_{I_j}$ for $1\le j\le k$, so that $M_k(\sigma)$ (as defined by \eqref{e:block.process}) is measurable with respect to $\filt_k$. Recall $\ell=(1-\gamma)L$, and define
	\[
	\Omega^{(\ell)}
	\equiv\bigg\{ \sigma\in\Omega :
		\dist\Big(M_\ell(\sigma),[a,b]\Big) \le
		L\gamma 
		\cdot \f{\RR\epsilon^{1/2}}{L^{1/2}}
		\cdot \NU^{1/2}
		\bigg\}\,.
	\]
Note that if $\sigma\notin\Omega^{(\ell)}$ then, by the definition 
\eqref{e:defn.RR} of $\RR$, we have
	\[
	|M_\ell(\sigma)|
	\ge L\gamma\cdot
		\f{\RR\epsilon^{1/2}}{L^{1/2}}\cdot \NU^{1/2}
		-\max\set{|a|,|b|}
	\ge L\gamma\cdot
		\f{\RR\epsilon^{1/2}}{2L^{1/2}}\cdot \NU^{1/2}\,,
	\]
where the last bound holds by using \eqref{e:defn.RR} together with the lower bound on $L$ from \eqref{e:restrict.L.general.disorder}, along with the assumption $\gamma=\epsilon$ and the fact that $\NU\ge1$ (see Assumption~\ref{a:subgaus}). It follows using Markov's inequality and the subgaussian tail bound that
	\begin{align*}
	\textup{(I)} &\equiv
	\P
	\bigg(
	\f{|\Omega^{(\ell)}|}{|\Omega|}
	\le \f12\bigg)
	\le \P
	\bigg(
	\f1{|\Omega|} \bigg|\bigg\{\sigma\in\Omega: 
		|M_\ell(\sigma)|\ge L\gamma
		\cdot\f{\RR\epsilon^{1/2}}{2L^{1/2}}\cdot \NU^{1/2}
		\bigg\}\bigg|
	\ge \f12
	\bigg)\\
	&\le 
	4\exp\bigg(- \f{L\RR^2 \gamma^2 \epsilon}{8} \bigg)
	\stackrel{\eqref{e:defn.RR}}{\le}
	4\exp\bigg( -\f{L\bar{r} \gamma\epsilon}{8} \bigg)
	\le \f{4}{\exp(L\epsilon^2)}\,,
	\end{align*}
where the last inequality again uses $\gamma=\epsilon$, and the fact that $\bar{r}\ge C'$ (a large absolute constant) by the definition \eqref{e:defn.RR}. Next, for $\ell+1\le k\le L$ let us define recursively the sets
	\[
	\Omega^{(k)}
	\equiv\bigg\{ \sigma\in\Omega^{(k-1)} :
		\dist\Big(M_k(\sigma),[a,b]\Big) \le
		(L-k)
		\cdot\f{\RR\epsilon^{1/2}}{L^{1/2}}
		\cdot\NU^{1/2}
		\bigg\}\,.
	\]
For each $\ell \le k \le L$ we further define the following bipartition of $\Omega^{(k)}$:
	\[\Omega^{(k),+}
	\equiv\bigg\{ \sigma\in\Omega^{(k)} :
		M_k(\sigma) \ge \f{a+b}{2}
		\bigg\}\,,\quad
	\Omega^{(k),-}
	\equiv \Omega^{(k)} \setminus \Omega^{(k),+}\,.\]
One of these subsets must be at least half the size of $\Omega^{(k)}$; without loss of generality we assume it is $\Omega^{(k),+}$. Assuming $|\Omega^{(k),+}|\ge \exp(K\delta)$, we can extract disjoint subsets $X_1,\ldots,X_a\subseteq\Omega^{(k),+}$ with $|X_a|=\exp(K\delta)$ such that the union of the $X_a$ occupies more than half the mass of $\Omega^{(k),+}$. Define
	\[
	\Xi^{(k)}
	\equiv \bigg\{\sigma\in\Omega^{(k),+}
	: \f{(\XI_{I_{k+1}},\sigma_{I_{k+1}})}{N^{1/2}}
	\in
	\bigg[\f{(\epsilon\NU)^{1/2} \RR}{L^{1/2}},\chi\bigg]
	\bigg\}\,.
	\]
It follows from the assumption that all pairs $\sigma\ne\tau$ in $\Omega$ are $\epsilon$-separated on $I_{k+1}$, so applying Lemma~\ref{l:Kbound.general} gives
	\beq\label{e:apply.lem.Kbound}
	\P\bigg(\f1m \sum_{a\le m} \I\bigg\{
	\f{|X_a\cap\Xi^{(k)}|}{|X_a|} 
	\le \exp\bigg( -\f{C''\NU s^2}{\epsilon}\bigg)
	\bigg\} \le \f12\bigg)
	\le \f{10}{\exp(s^2\epsilon)}\,,\eeq
provided that the parameters satisfy the bounds
	\beq\label{e:interval.gaus.restrictions}
	\RR \le s \le 
	s_{\max} \equiv
	\f{(L\epsilon)^{1/2}\chi}{2(C'')^{1/2}\NU}
	\le s_2
	\equiv\f1{C''} \bigg( \f{K\epsilon^3}{\NU} \bigg)^{1/4}
	\,.\eeq
Note that in the restriction \eqref{e:restrict.L.general.disorder},
the lower bound on $L$ guarantees that we in fact have $s_{\max}\ge\bar{R}$ above, while the upper bound on $L$ guarantees $s_{\max}\le s_2$. Note that $X_a \cap \Xi^{(k)}\subseteq \Omega^{(k+1)}$ for all $a\le m$. It follows that, on the event $|\Omega^{(k)}|\ge 2\exp(K\delta)$, we have
	\beq\label{e:interval.tail.bound.general.disorder}
	\P\bigg(
	\f{|\Omega^{(k+1)}|}{|\Omega^{(k)}|}
	\le \f14\exp\bigg( -\f{C''\NU s^2}{\epsilon}\bigg)
	\,\bigg|\,\filt_k 
	\bigg)
	\le \f{10}{\exp(s^2\epsilon)}\,.\eeq
To apply this bound, we define the stopping time $\ell\le\ZETA\le L$ by
	\beq\label{e:stopping.time.zeta.general.disorder}
	\ZETA
	\equiv
	\min\bigg\{ k\ge\ell : 
	k=L, \textup{ or }
	\f{|\Omega^{(k)}|}{2}
	<\exp(K\delta) 
	\bigg\}\,.\eeq
Applying \eqref{e:interval.tail.bound.general.disorder} with $s=s_{\max}$ and $\ell\le k<\ZETA$ gives, on the event $|\Omega^{(\ell)}|\ge|\Omega|/2$, the bound
	\begin{align*}
	\textup{(II)}
	&\equiv\P\bigg(
	\min_{\ell\le k<\ZETA}
	\f{|\Omega^{(k+1)}|}{|\Omega^{(k)}|}
	\le \f14 \exp\bigg( -\f{C''\NU (s_{\max})^2}{\epsilon}\bigg)
	\,\bigg|\,\filt_\ell
	\bigg) \\
	&\stackrel{\eqref{e:interval.tail.bound.general.disorder}}{\le}
	\f{L\gamma \cdot 10}{\exp((s_{\max})^2\epsilon)}
	=10\gamma\exp\bigg( \log L- 
		\f{L\epsilon\chi^2}{4C''\NU^2} \bigg)
	\le
	10\gamma\exp\bigg( -
		\f{L\epsilon\chi^2}{8C''\NU^2} \bigg)
	\,,
	\end{align*}
where the last bound again uses the lower bound on $L$ from \eqref{e:restrict.L.general.disorder}.
We then define the truncated random variables
	\[
	V_{k+1}\equiv \min
	\bigg\{\f{|\Omega^{(k)}|}{|\Omega^{(k+1)}|},
	 4\exp\bigg( \f{C''\NU (s_{\max})^2}{\epsilon}\bigg)
	\bigg\}\,.
	\]
It follows by simplifying \eqref{e:interval.tail.bound.general.disorder} that, on the event $|\Omega^{(k)}|\ge2\exp(K\delta)$, we have for all $s\ge \RR$ that
	\[
	\P\bigg(V_{k+1} \ge
		 \exp\bigg( \f{2C''\NU s^2}{\epsilon}\bigg)
		\,\bigg|\,\filt_k\bigg)
	\le \f{10}{\exp(s^2\epsilon)}\,.
	\]
(More precisely, for $\RR \le s\le s_{\max}$, the above bound holds due to \eqref{e:interval.tail.bound.general.disorder}. For $s\ge s_{\max}$, the bound holds trivially, because the truncation in the definition of $V_{k+1}$ implies that the left-hand side above is zero.) Making the change of variables $t=\exp(s^2\epsilon/2)$ gives, for all $t\ge t_{\min}=\exp(\RR^2\epsilon/2)$, 
	\[
	\P\bigg( (V_{k+1})^{\epsilon^2/(4C''\NU)} \ge t
		\,\bigg|\, \filt_k\bigg)
	\le \f{10}{t^2}\,.
	\]
Recall from the definition \eqref{e:defn.RR}, and the assumption $\gamma=\epsilon$, that $\RR^2\epsilon=\bar{r}\ge C'$ (a large absolute constant). Integrating the above gives the moment bound
	\beq\label{e:sqrt.V.mmt.bound.general.disorder}
	\E\bigg[ (V_{k+1})^{\epsilon^2/(4C''\NU)} \,\bigg|\, \filt_k\bigg]
	\le t_{\min} + 
		\int_{t_{\min}}^\infty \f{10}{t^2}\,dt
	= t_{\min} + \f{10}{t_{\min}}
	\le \exp(\RR^2\epsilon)\,,\eeq
where the last inequality uses that $\RR^2\epsilon=\bar{r}$ is large by \eqref{e:defn.RR}. Consequently, applying Markov's inequality we obtain, on the event $|\Omega^{(\ell)}|\ge|\Omega|/2$, for any $w\ge0$,
	\[\P\bigg(
	\prod_{\ell<k\le\ZETA} V_k 
	\ge 
	\exp\bigg( \f{4 C''\NU w^2}{\epsilon} \bigg)
	\,\bigg|\,\filt_\ell \bigg)
	\stackrel{\eqref{e:sqrt.V.mmt.bound.general.disorder}}{\le} 
	\f{\exp(L\gamma \cdot \RR^2 \epsilon   )
	}{\exp(w^2\epsilon)}\,.\]
Setting $w^2= 2L\RR^2 \gamma$ then gives the bound 
	\[
	\textup{(III)}
	\equiv
	\P\bigg(
	\prod_{\ell<k\le\ZETA} V_k 
	\ge 
	\exp\bigg( \f{8 C''\NU L\RR^2 \gamma}{\epsilon} \bigg)
	\,\bigg|\,\filt_\ell \bigg)
	\le
	\f1{\exp(L\RR^2\gamma\epsilon)}
	\stackrel{\eqref{e:defn.RR}}{=}
	\f1{\exp(L\bar{r}\epsilon)}
	\,.\]
Finally, let us note that on the event
	\[\bm{E}=\bigg\{
	\f{|\Omega^{(\ell)}|}{|\Omega|}
	\ge \f12\,,\min_{\ell\le k<\ZETA}
	\f{|\Omega^{(k+1)}|}{|\Omega^{(k)}|}
	\ge \f14
	\exp\bigg(- \f{C''\NU (s_{\max})^2}{\epsilon}
	\bigg),
	\prod_{\ell<k\le\ZETA} V_k 
	\le\exp\bigg( \f{8C''\NU L\RR^2 \gamma}{\epsilon}\bigg)
	\bigg\}\,,
	\]
we must have $|\Omega^{(k+1)}|/|\Omega^{(k)}|=1/V_{k+1}$ for all $\ell\le k<\ZETA$. Consequently, recalling $K=N/L$, we have 
	\[
	\f{|\Omega^{(\ZETA)}|}{2}
	= \f{|\Omega|}{2}
	\cdot\f{|\Omega^{(\ell)}|}{|\Omega|}
	\cdot\prod_{\ell<k\le\ZETA} \f1{V_k}
	\ge \f14\exp\bigg( 2K\delta
		- \f{8C''\NU L\RR^2 \gamma}{\epsilon}\bigg)
	\ge \exp(K\delta)\,,
	\]
where the last inequality uses the upper bound on $L$ from \eqref{e:restrict.L.general.disorder}.
From the definition \eqref{e:stopping.time.zeta.general.disorder} of $\ZETA$, this implies $\ZETA=L$. Combining the above bounds gives
\begin{align*}
	\textup{(IV)}&\equiv\P\bigg(
	\f{|\Omega^{(L)}|}{|\Omega|}
	\le \f14
	\exp\bigg(-\f{8C''\NU L\RR^2 \gamma}{\epsilon}\bigg)
	\bigg) \le 1-\P(\bm{E}) 
	\le \textup{(I)}+\textup{(II)}+\textup{(III)}\\
	&\le \f{4}{\exp(L\epsilon^2)}
	+10\gamma\exp\bigg( - \f{L\epsilon\chi^2}{8C''\NU^2} \bigg)
	+\f1{\exp(L\bar{r}\epsilon)}
	\le \exp\bigg(-\f{L\epsilon^2}{2}\bigg)\,,
	\end{align*}
where the last bound holds for $\epsilon$ small enough (which is guaranteed by taking $\delta$ small enough). Finally, recall that $\sigma\in\Omega^{(L)}$ implies $M_L(\sigma)=(\XI,\sigma)/N^{1/2}\in[a,b]$ as desired, so this concludes the proof. 
\end{proof}
\end{ppn}

\begin{proof}[\hypertarget{proof:t.interval.general}{Proof of Theorem~\ref{t:interval.general}}] Since we can adjust the constants later, we can assume $|S|\ge\exp(9N\delta)$. Recalling the statement of Proposition~\ref{p:sep.on.blocks.subsample}, let $\epsilon,\gamma$ be positive constants with $\psi_2(8\epsilon)=\delta$ and $\gamma=\epsilon$. Assume  $L$ satisfies \eqref{e:restrict.L.general.disorder}. It follows from Proposition~\ref{p:sep.on.blocks.subsample} that there exists $J_1\subseteq[L]$ with $|J_1|=L\gamma$, and $\Omega_1\subseteq S$ with $|\Omega_1|=\exp(N\delta/L)$, such that all pairs $\sigma\ne\tau$ in $\Omega_1$ are $\epsilon$-separated on $I_j$ for all $j\in J_1$. We can apply Proposition~\ref{p:sep.on.blocks.subsample} again on $S\setminus\Omega_1$, and so on, to extract
disjoint subsets $\Omega_1,\ldots,\Omega_m\subseteq S$, up to the first $m$ such that
	\[
	\sum_{a\le m}|\Omega_a| \ge \f{|S|}{2}\,.
	\]
Then for each $a$ we will have $|\Omega_a|=\exp(N\delta/L)$,
and all pairs $\sigma\ne\tau$ in $\Omega_a$ will be $\epsilon$-separated on $I_j$ for all $j\in J_a$, where $J_a\subseteq[L]$ with $|J_a|=L\gamma$. Let
	\[
	\Xi\equiv\bigg\{\sigma\in S: 
		\f{(g,\sigma)}{N^{1/2}}\in[a,b]\bigg\}\,,\quad
	\Xi_a
	\equiv \Xi \cap \Omega_a\,.
	\]
By Markov's inequality and Proposition~\ref{p:addinterval.general},
	\[
	\P\bigg(\f1m \sum_{a\le m}
	\I\bigg\{ \f{|\Xi_a|}{|\Omega_a|}
	\le\f14\exp\bigg( -\f{8C''\NU \bar{r} L}{\epsilon}\bigg) 
	\bigg\} \ge \f12\bigg)
	\le \exp\bigg(-\f{L\epsilon^2}{2}\bigg)\,.
	\]
On the complementary event we must have
	\[
	|\Xi|
	\ge \sum_{a\le m} |\Xi_a|
	\ge \f{m|\Omega_1|}{2}
	\cdot\f14\exp\bigg( -\f{8C''\NU \bar{r} L}{\epsilon}\bigg) 
	\ge \f{|S|}{16}
	\exp\bigg( -\f{8C''\NU \bar{r} L}{\epsilon}\bigg) 
	\,,
	\]
so we have shown that, for any $L$ satisfying \eqref{e:restrict.L.general.disorder}, we have
	\[
	\P\bigg(
	\f{|\Xi|}{|S|}
	\le \f1{16}\exp\bigg( -\f{8C''\bar{r}\NU L}{\epsilon}\bigg) 
	\bigg)
	\le \exp\bigg( -\f{L\epsilon^2}{2}\bigg)\,.
	\]
The conclusion follows by recalling that we chose $\gamma$ and $\epsilon$ depending on $\delta$. 
\end{proof}

\begin{proof}[{\hypertarget{proof:t.interval.gaus}{Proof of Theorem~\ref{t:interval.gaus}}}] This follows essentially by the \hyperlink{proof:t.interval.general}{proof of Theorem~\ref{t:interval.general}}, the main difference being that in the proof of Proposition~\ref{p:addinterval.general} we can apply Lemma~\ref{l:Kbound.gaus} in place of Lemma~\ref{l:Kbound.general}. Thus in place of \eqref{e:apply.lem.Kbound} we will have
	\[\P\bigg(\f1m \sum_{a\le m} \I\bigg\{
	\f{|X_a\cap\Xi^{(k)}|}{|X_a|} 
	\le \f1{\exp(2C's^2)}
	\bigg\} \le \f12\bigg)
	\le \f{4}{\exp(s^2\epsilon/C')}\,,\]
where instead of \eqref{e:interval.gaus.restrictions} we will require (from the conditions of Lemma~\ref{l:Kbound.gaus})
	\[\RR \le s \le 
	s_{\max}
	\equiv
	\f{L^{1/2}\chi}{(6C')^{1/2}}
	\le s_2
	\equiv \f{(K\delta)^{1/2}}{C'}\,.
	\]
Since $KL=N$, this explains why we ultimately require $w\le N^{1/2}/C_\delta$, in contrast with Theorem~\ref{t:interval.general} which requires $w\le N^{1/3}/C_\delta$.
\end{proof}

\section{Universality on average in perceptron models}
\label{s:univ}

Most of this section is devoted to the proof of the following theorem, which says that the \textbf{expectation} of the perceptron free energy (with a suitable truncation) is universal with respect to the disorder. In Section~\ref{s:sharp} we will combine this with a concentration result (Proposition~\ref{p:truncated.conc}) to yield the \hyperlink{proof:t.univ}{proof of Theorem~\ref{t:univ}}. 

\begin{thm}\label{t:univ.avg}
Consider a perceptron model \eqref{e:Z.U} where $U$ is $\set{0,1}$-valued and piecewise continuous. Let $Z(g)$ be the partition function with standard gaussian disorder, and let $Z(\XI)$ be the partition function with general disorder, assuming the $(\XI^a)_i$ are i.i.d.\ random variables with zero mean, unit variance, and finite third moment. Then
	\[
	\f1N \Big| \E\log_{N\delta}Z(g)
	-\E\log_{N\delta}Z(\XI)\Big| = o_N(1)\,.
	\]
for any small positive constant $\delta$.
\end{thm}

This section is organized as follows:
\begin{itemize}
\item In \S\ref{ss:clt} we give some multivariate central limit theorem estimates. The main estimate  is Corollary~\ref{c:clt}.

\item In \S\ref{ss:univ.postemp} we use Corollary~\ref{c:clt} to prove a version (Proposition~\ref{p:univ.positive.temp}) of Theorem~\ref{t:univ.avg} for a ``positive-temperature'' variant of the model.

\item In \S\ref{ss:postemp.zerotemp} we show how to transfer 
the result of Proposition~\ref{p:univ.positive.temp} from positive temperature to zero temperature, leading to the 
\hyperlink{proof:t.univ.avg}{proof of Theorem~\ref{t:univ.avg}}.
\end{itemize}

\subsection{Central limit theorem estimates}
\label{ss:clt}

In this subsection we state and prove a consequence of the multivariate central limit theorem,
Corollary~\ref{c:clt} below,
which was already used in 
the proof of Proposition~\ref{p:general.disorder.threshold}.
Corollary~\ref{c:clt} will also be used below in the proof of Proposition~\ref{p:univ.positive.temp}, which is a preliminary universality result. We begin with a central limit theorem for sufficiently smooth functions, which is similar to \cite[Thm.~4.1]{MR2594620}: 

\begin{lem}\label{l:clt.smooth}
Let $g$ be a standard gaussian vector in $\R^N$. Let $\XI$ be a random vector in $\R^N$ with i.i.d.\ coordinates $\XI_i$, with mean zero, unit variance, and finite third moment. If $F:\R^p\to\R$ is a bounded function with derivatives up to third order bounded uniformly by a finite constant $C_F$, then we have
	\[
	\bigg|\E\bigg[ 
	F\bigg( \f{\Sigma\XI}{N^{1/2}}\bigg)
	-F\bigg( \f{\Sigma g}{N^{1/2}}\bigg)
	\bigg]\bigg|
	\le \f{pC_F}{N^{1/2}}
	\bigg\{\E(|\XI_1|^3)+\E(|g_1|^3)\bigg\}
	\]
for any matrix $\Sigma\in\set{-1,+1}^{p\times N}$.

\begin{proof}
Write $\sigma^{(\ell)}\in\set{-1,+1}^N$ for the $\ell$-th row of $\Sigma$.
Let $\Sigma_{j,\circ}$ denote the $p\times N$ matrix that results from 
zeroing the $j$-th row of $\Sigma$. We interpolate between $g$ and $\XI$ by defining  $\bXI^j\equiv(\XI_1,\ldots,\XI_j,g_{j+1},\ldots,g_N)$. Let $\bXI^{j,\circ}$ denote the vector that results from zeroing the $j$-th entry of $\bXI^j$. We then define the $p$-dimensional vectors
	\begin{align*}
	\bDe_j\equiv \f{\Sigma\bXI^j}{N^{1/2}}
	&= \bigg( 
	\sum_{i=1}^j \f{\XI_i (\sigma^{(\ell)})_i}{N^{1/2}}
	+\sum_{i=j+1}^N \f{g_i (\sigma^{(\ell)})_i}{N^{1/2}}
	\bigg)_{\ell\le p}\,, \\
	\bDe_{j,\circ}\equiv \f{\Sigma^{j,\circ} \bXI^{j,\circ}}{N^{1/2}}
	&= \bigg( 
	\sum_{i=1}^{j-1} \f{\XI_i (\sigma^{(\ell)})_i}{N^{1/2}}
	+\sum_{i=j+1}^N \f{g_i (\sigma^{(\ell)})_i}{N^{1/2}}
	\bigg)_{\ell\le p}\,.\end{align*}
Then the quantity of interest can be bounded as
	\[\bigg|\E\Big[ F(\bDe_N)-F(\bDe_0)\Big]\bigg|
	\le\sum_{j\le N}
	\bigg| \E\Big[ F(\bDe_j)-F(\bDe_{j-1})\Big]\bigg|\,,
	\]
so it suffices to bound $|\E[ F(\bDe_j)-F(\bDe_{j-1})]|$ for each $1\le j\le N$. To this end we let
	\begin{align*}
	\bDe_{j,\ell}
	&\equiv 
	\bDe_{j,\circ}+\f1{N^{1/2}}
	\bigg( \XI_j(\sigma^{(1)})_j,\ldots,\XI_j(\sigma^{(\ell)})_j,
	g_j(\sigma^{(\ell+1)})_j,\ldots,g_j(\sigma^{(p)})_j
	\bigg)\,,\\
	\bDe_{j,\ell,\circ}
	&\equiv
	\bDe_{j,\circ}+\f1{N^{1/2}}
	\bigg( \XI_j(\sigma^{(1)})_j,\ldots,\XI_j(\sigma^{(\ell-1)})_j,
	0,
	g_j(\sigma^{(\ell+1)})_j,\ldots,g_j(\sigma^{(p)})_j
	\bigg)\,,
	\end{align*}
and note that $F(\bDe_{j-1})=F(\bDe_{j,0})$ while $F(\bDe_j)=F(\bDe_{j,p})$. It follows that
	\[
	\bigg|\E\Big[ F(\bDe_j)-F(\bDe_{j-1})\Big]\bigg| 
	\le \sum_{\ell=1}^p
	\bigg| \E\Big[ 
	F(\bDe_{j,\ell})-F(\bDe_{j,\ell-1}) \Big]\bigg|\,.
	\]
We then apply the assumption on $F$ to Taylor expand
	\begin{align*}
	F(\bDe_{j,\ell})
	&= F(\bDe_{j,\ell,\circ})
	+ \partial_\ell F(\bDe_{j,\ell,\circ})
	\f{\XI_j (\sigma^{(\ell)})_j}{N^{1/2}}
	+ \f{(\partial_\ell)^2 F(\bDe_{j,\ell,\circ})}{2}
	\f{(\XI_j)^2}{N}
	+\f{R^\XI}{N^{3/2}}\,,\\
	F(\bDe_{j,\ell-1})
	&= F(\bDe_{j,\ell,\circ})
	+\partial_\ell F(\bDe_{j,\ell,\circ})
	\f{g_j (\sigma^{(\ell)})_j}{N^{1/2}}
	+\f{(\partial_\ell)^2 F(\bDe_{j,\ell,\circ})}2
	\f{(g_j)^2}{N}+ \f{R^g}{N^{3/2}}\,,
	\end{align*}
where $\E|R^\XI|\le C_F\E(|\XI_j|^3)$ and $\E|R^g|\le C_F\E(|g_j|^3)$. Matching the first and second moments of $\XI_j$ and $g_j$ gives
	\[
	\bigg| \E\Big[ 
	F(\bDe_{j,\ell})-F(\bDe_{j,\ell-1}) \Big]\bigg|
	\le \f{C_F}{N^{3/2}}
	\bigg\{\E\Big(|\XI_j|^3+|g_j|^3\Big)\bigg\}\,.\]
Summing over $1\le\ell\le p$ and $1\le j\le N$ gives the desired overall bound on 
$|\E[F(\bDe_N)-F(\bDe_0)]|$.
\end{proof}
\end{lem}

\begin{lem}\label{l:clt.cts}
Let $g$ be a standard gaussian vector in $\R^N$. Let $\XI$ be a random vector in $\R^N$ with i.i.d.\ coordinates $\XI_i$, with mean zero, unit variance, and finite third moment. If $F:\R^p\to\R$ is a bounded  continuous function, then
	\[\adjustlimits
	\lim_{N\to\infty}
	\max\bigg\{
	\bigg|\E\bigg[ 
	F\bigg( \f{\Sigma\XI}{N^{1/2}}\bigg)
	-F\bigg( \f{\Sigma g}{N^{1/2}}\bigg)
	\bigg]\bigg|
	: \Sigma\in\set{-1,+1}^{p\times N}
	\bigg\}=0\,.
	\]

\begin{proof}
For $\eta>0$ define the smoothed function
$F_\eta(x)\equiv \E F(x+\eta \bm{z})$
where $\bm{z}$ is a standard gaussian vector in $\R^p$. For any $\eta>0$ the function $F_\eta$ satisfies the conditions of Lemma~\ref{l:clt.smooth}, so we have
	\beq\label{e:smoothed.clt}
	\textup{(I)}
	\equiv \max\bigg\{
	\bigg|\E\bigg[ 
	F_\eta\bigg( \f{\Sigma\XI}{N^{1/2}}\bigg)
	-F_\eta\bigg( \f{\Sigma g}{N^{1/2}}\bigg)
	\bigg]\bigg|
	: \Sigma\in\set{-1,+1}^N
	\bigg\}
	\le \f{p C_{F,\eta}}{N^{1/2}}
	\bigg\{\E\Big(|\XI_j|^3+|g_j|^3\Big)\bigg\}\,.
	\eeq
Since $F$ is bounded and continuous, the function $F_\eta$ converges locally uniformly to $F$ as $\eta\downarrow0$, so we have
	\begin{align*}
	\textup{(II)}
	&\equiv
	\max\bigg\{
	\E\bigg[ \bigg|
	F_\eta\bigg( \f{\Sigma\XI}{N^{1/2}}\bigg)
	-F\bigg( \f{\Sigma\XI}{N^{1/2}}\bigg)
	\bigg|
	; \f{\|\Sigma\XI\|_\infty}{N^{1/2}} \le R\bigg]
	: \Sigma\in\set{-1,+1}^N
	\bigg\}\\
	&\le 
	\sup\bigg\{ |F_\eta(x)-F(x)|
	: \|x\|_\infty\le R\bigg\}
	\stackrel{\eta\downarrow0}{\longrightarrow} 0
	\end{align*}
for any finite $R$. For the complementary event $\|\Sigma\XI\|_\infty> R$, we can use Chebychev's inequality to bound
	\begin{align*}
	\textup{(III)}
	&\equiv\max\bigg\{
	\E\bigg[ \bigg|
	F_\eta\bigg( \f{\Sigma\XI}{N^{1/2}}\bigg)
	-F\bigg( \f{\Sigma\XI}{N^{1/2}}\bigg)
	\bigg|
	; \f{\|\Sigma\XI\|_\infty}{N^{1/2}} > R\bigg]
	: \Sigma\in\set{-1,+1}^N
	\bigg\} \\
	&\le 2\|F\|_\infty
	\P\bigg(\f{\|\Sigma\XI\|_\infty}{N^{1/2}} > R\bigg)
	\le \f{2p\|F\|_\infty }{R^2}\,,
	\end{align*}
which can be made arbitrarily small by choosing $R=R(p)$ large. The claim follows by combining the estimates for $\textup{(I)}$, $\textup{(II)}$, and $\textup{(III)}$.\end{proof}
\end{lem}

\begin{cor}\label{c:clt} Let $g$ be a standard gaussian vector in $\R^N$. Let $\XI$ be a random vector in $\R^N$ with i.i.d.\ coordinates $\XI_i$, with mean zero, unit variance, and finite third moment. 
If $f:\R\to\R$ is a bounded piecewise continuous function, and we define $F:\R^p\to\R$ by $F(x_1,\ldots,x_p)=f(x_1)\cdots f(x_p)$, then
	\[
	\lim_{N\to\infty}
	\max\bigg\{
	\bigg|\E\bigg[ 
	F\bigg( \f{\Sigma\XI}{N^{1/2}}\bigg)
	-F\bigg( \f{\Sigma g}{N^{1/2}}\bigg)
	\bigg]\bigg|
	: \Sigma\in\set{-1,+1}^{p\times N}
	\bigg\}
	=0\,.
	\]

\begin{proof}
As in the proof of Lemma~\ref{l:clt.cts}, consider the smoothed function
	\[
	F_\eta(x)
	=\E F(x+\eta\bm{z})
	= \E\bigg( \prod_{\ell=1}^p f(x_\ell + \eta z_\ell)\bigg)
	=  \prod_{\ell=1}^p f_\eta(x_\ell)\,.
	\]
We then have the approximation result of Lemma~\ref{l:clt.smooth}, which gives (as in \eqref{e:smoothed.clt})
	\beq\label{e:smoothed.clt.REPEAT}
	\max\bigg\{
	\bigg|\E\bigg[ 
	F_\eta\bigg( \f{\Sigma\XI}{N^{1/2}}\bigg)
	-F_\eta\bigg( \f{\Sigma g}{N^{1/2}}\bigg)
	\bigg]\bigg|
	: \Sigma\in\set{-1,+1}^N
	\bigg\}
	\le \f{p C_{F,\eta}}{N^{1/2}}
	\bigg\{\E\Big(|\XI_j|^3+|g_j|^3\Big)\bigg\}\,.\eeq
We next bound
	\beq\label{e:clt.product}
	\E\bigg[ \bigg|F_\eta\bigg(\f{\Sigma \XI}{N^{1/2}}\bigg)
	-F\bigg(\f{\Sigma \XI}{N^{1/2}}\bigg)\bigg|
	\bigg]
	\le  (\|f\|_\infty)^{p-1}\sum_{\ell=1}^p
	\E\bigg[\bigg|
	f_\eta\bigg( \f{(\XI,\sigma^{(\ell)})}{N^{1/2}}\bigg)
	-f\bigg( \f{(\XI,\sigma^{(\ell)})}{N^{1/2}}\bigg) \bigg|\bigg]\,,
	\eeq
and likewise with $g$ in place of $\XI$. Recall the assumption that $f$ is piecewise continuous, so its set of discontinuities $D_f$ is finite.
Therefore, given any $\theta>0$, we can construct a continuous function $\Upsilon:\R\to[0,1]$ such that $\Upsilon=1$ on an open set $O_f\supseteq D_f$ with $\R\setminus O_f$ compact, but $\E \Upsilon(g_1)\le\theta$ if $g_1$ is a standard gaussian random variable. Applying Lemma~\ref{l:clt.cts} gives, for any fixed choice of $O_f$ and $\Upsilon$,
	\[
	\Delta_N\equiv
	\max\bigg\{
	\bigg|\E \bigg[\Upsilon\bigg( \f{(\XI,\sigma)}{N^{1/2}}\bigg)
	- \Upsilon\bigg( \f{(g,\sigma)}{N^{1/2}}\bigg)
	\bigg]\bigg| : \sigma\in\set{-1,+1}^N
	\bigg\} \stackrel{N\to\infty}{\longrightarrow} 0\,.
	\]
For any fixed $\sigma\in\set{-1,+1}^N$, the scalar product $(g,\sigma)/N^{1/2}$ is a standard gaussian random variable. It follows that
	\begin{align*}
	\textup{(I)}
	&\equiv\max\bigg\{\E\bigg[\bigg|
	f_\eta\bigg( \f{(\XI,\sigma^{(\ell)})}{N^{1/2}}\bigg)
	-f\bigg( \f{(\XI,\sigma^{(\ell)})}{N^{1/2}}\bigg) \bigg|
	; \f{(\XI,\sigma^{(\ell)})}{N^{1/2}} \in O_f\bigg]
	: \sigma^{(\ell)} \in \set{-1,+1}^N
	\bigg\}\\
	&\le 2\|f\|_\infty
	\E \bigg[\Upsilon\bigg( \f{(\XI,\sigma)}{N^{1/2}}\bigg)\bigg]
	\le 2\|f\|_\infty
	\bigg\{ \Delta_N
	+\E \bigg[\Upsilon\bigg( \f{(g,\sigma)}{N^{1/2}}\bigg)\bigg]
	\bigg\} 
	\le  2\|f\|_\infty
	\bigg\{ \Delta_N + \theta\bigg\}
	\end{align*}
We have $f_\eta$ converging to $f$ uniformly on the compact set $\R\setminus O_f$, so
	\[\textup{(II)}
	\equiv\max\bigg\{\E\bigg[\bigg|
	f_\eta\bigg( \f{(\XI,\sigma^{(\ell)})}{N^{1/2}}\bigg)
	-f\bigg( \f{(\XI,\sigma^{(\ell)})}{N^{1/2}}\bigg) \bigg|
	; \f{(\XI,\sigma^{(\ell)})}{N^{1/2}} \in \R
	\setminus O_f\bigg]
	: \sigma^{(\ell)}\in \set{-1,+1}^N
	\bigg\}
	\stackrel{\eta\downarrow0}{\longrightarrow} 0\,.
	\]
Of course, the estimates for \textup{(I)} and \textup{(II)} also hold with $g$ in place of $\XI$.  The claim follows by substituting those estimates into \eqref{e:clt.product}, and combining with \eqref{e:smoothed.clt.REPEAT}.\end{proof}
\end{cor}

\subsection{Universality at positive temperature} 
\label{ss:univ.postemp}

In this subsection we 
establish universality for a perceptron model at ``positive temperature,'' more precisely, of the form \eqref{e:Z.U} where $\log U:\R\to (-\infty,0]$ is uniformly bounded. Note that the Azuma--Hoeffding martingale inequality implies in this case that $N^{-1}\log Z(U;\XI)$ is exponentially well concentrated around its mean, although we will not ultimately use this fact. The following result shows that the mean value is universal with respect to the disorder model:

\begin{ppn}\label{p:univ.positive.temp}
For the ``positive-temperature'' perceptron model \eqref{e:Z.U} where $\log U:\R\to (-\infty,0]$ is uniformly bounded and piecewise continuous, let $Z(U;g)$ be the partition function with standard gaussian disorder,
and let $Z(U;\XI)$ be the partition function with general disorder, assuming the $(\XI^a)_i$ are i.i.d.\ random variables with zero mean, unit variance, and finite third moment. Then
	\[
	\lim_{N\to\infty}\f1N \Big| \E\log Z(U;g)
	-\E\log Z(U;\XI)\Big|=0\,.
	\]

\begin{proof} It follows from Corollary~\ref{c:clt} that, uniformly over $\sigma^{(1)},\ldots,\sigma^{(\ell)}\in\set{-1,+1}^N$, we have 
	\beq\label{e:apply.clt.prod}
	\Delta_j(\sigma^{(1)},\ldots,\sigma^{(\ell)})
	\equiv\E\bigg[
	\prod_{k=1}^\ell U \bigg( \f{(\XI^j,\sigma^{(k)} )}{N^{1/2}}\bigg)
	-\prod_{k=1}^\ell U \bigg( \f{(g^j,\sigma^{(k)} )}{N^{1/2}}\bigg)
	\bigg]
	\stackrel{N\to\infty}{\longrightarrow} 0\,.
	\eeq
We interpolate between $Z(U;g)$ and $Z(U;\XI)$ by defining
	\begin{align*}
	Z_j &\equiv \sum_{\sigma\in\set{-1,+1}^N}
	\prod_{i=1}^j
		U\bigg( \f{(\XI^i,\sigma)}{N^{1/2}}\bigg)
	\prod_{i=j+1}^M
	U\bigg( \f{(g^i,\sigma)}{N^{1/2}}\bigg)\bigg\}\,,\\
	Z_{j,\circ}
	&\equiv\sum_{\sigma\in\set{-1,+1}^N}
	\prod_{i=1}^{j-1} U\bigg( \f{(\XI^i,\sigma)}{N^{1/2}}\bigg)
	\prod_{i=j+1}^M U\bigg( \f{(g^i,\sigma)}{N^{1/2}}\bigg)\bigg\}\,.
	\end{align*}
Then $Z_0=Z(U;g)$ and $Z_M=Z(U;\XI)$, so
	\[
	\log Z(U;\XI)-\log Z(U;g)
	= \sum_{j\le M} \Big( \log Z_j -\log Z_{j-1}\Big)
	\equiv \sum_{j\le M} Y_j\,.
	\]
Writing $\<\cdot\>_{j,\circ}$ for expectation with respect to the Gibbs measure $\mu_{j,\circ}$ corresponding to $Z_{j,\circ}$, we have
	\[
	Y_j
	= \log \f{Z_j}{Z_{j,\circ}}
		-\log \f{Z_{j-1}}{Z_{j,\circ}}
	=\log \bigg\< U \bigg( \f{(\XI^j,\sigma)}{N^{1/2}}\bigg)\bigg\>_{j,\circ}
	-\log \bigg\< U \bigg( \f{(g^j,\sigma)}{N^{1/2}}\bigg)\bigg\>_{j,\circ}
	\]
Given any $A<\infty$ and $\eta>0$ we can choose $p=p(A,\eta)<\infty$ such that
	\[
	\sup\bigg\{ \bigg| \log(1+x)
	- \sum_{\ell\le p} \f{(-1)^{\ell+1}}{\ell} x^\ell\bigg|
	: \f1{\exp(A)} \le 1+ x\le 0 \bigg\} \le\eta
	\]
It follows that $|\E Y_\ell - \bar{y}_\ell|\le \eta$ where
	\[
	\bar{y}_j
	\equiv
	\sum_{\ell\le p}
	\f{(-1)^{\ell+1}}{\ell}
	\E\bigg[
	\bigg(\bigg\< U \bigg( \f{(\XI^j,\sigma)}{N^{1/2}}\bigg)\bigg\>_{j,\circ}-1\bigg)^\ell
	-\bigg(\bigg\< U \bigg( \f{(g^j,\sigma)}{N^{1/2}}\bigg)\bigg\>_{j,\circ}-1\bigg)^\ell\bigg]\,.
	\]
Let $c_{p,\ell}$ denote the coefficients such that
	\[
	\sum_{\ell\le p}
	\f{(-1)^{\ell+1}}{\ell}(x-1)^\ell
	= \sum_{\ell\le p} c_{p,\ell} x^\ell\,.
	\]
Recalling that $\mu_{j,\circ}$ is the Gibbs measure corresponding to $Z_{j,\circ}$, we can rewrite
	\begin{align*}
	\bar{y}_j
	&= \sum_{\ell\le p}
	c_{p,\ell}
	\E\bigg[ \bigg(\bigg\< U \bigg( \f{(\XI^j,\sigma)}{N^{1/2}}\bigg)\bigg\>_{j,\circ}\bigg)^\ell
	-\bigg(\bigg\< U \bigg( \f{(g^j,\sigma)}{N^{1/2}}\bigg)\bigg\>_{j,\circ}\bigg)^\ell
	\bigg] \\
	&=  \sum_{\ell\le p}
	c_{p,\ell}
	\sum_{\sigma^{(1)},\ldots,\sigma^{(\ell)}}
	\Delta_j(\sigma^{(1)},\ldots,\sigma^{(\ell)})
	\prod_{k=1}^\ell \mu_{j,\circ}(\sigma^{(k)})\,,
	\end{align*}
which by \eqref{e:apply.clt.prod} converges to zero. The claim follows by recalling that $|\E Y_\ell - \bar{y}_\ell|\le \eta$, and $\eta>0$ was arbitrary. \end{proof}
\end{ppn}

\subsection{Positive temperature to zero temperature} 
\label{ss:postemp.zerotemp}

We first review a simple inequality \eqref{e:log.delta.bound} which was previously proved by \cite[Lem.~8.3.7]{MR3024566}. If $x\ge y\ge 0$, then for any $\Gamma<0$ we have
	\[
	\log_\Gamma y-\log_\Gamma x
	= \begin{cases}
	\log y-\log x & \textup{if $x \ge y \ge \exp(\Gamma)$,}\\
	\Gamma-\log x & \textup{if $x\ge \exp(\Gamma)\ge y$,}\\
	0 & \textup{if $\exp(\Gamma)\ge x\ge y$}.
	\end{cases}
	\]
If we further assume $1\ge x\ge y\ge0$, then the above quantities are all in $[\Gamma,0]$. It follows that 
	\beq\label{e:log.delta.bound}
	0\ge \log_\Gamma y-\log_\Gamma x
	\ge \max\bigg\{ \log\f{y}{x},\Gamma\bigg\}
	=\log_\Gamma\f{y}{x}\,.
	\eeq
We will use this inequality in the proof of the next result.

\begin{ppn}\label{p:pos.to.zero.temp}
Consider a perceptron model \eqref{e:Z.U} where $U$ is $\set{0,1}$-valued and piecewise continuous. Let $Z(\XI)$ denote the partition function with general disorder, assuming the $(\XI^a)_i$ are i.i.d.\ random variables with zero mean, unit variance, and finite third moment. Let $u_A(x) \equiv \max\set{\log U(x), -A}$. Then for any positive constant $\delta$ we have
	\[
	\f1N\E\Big| \log_{N\delta} Z(u_A;\XI) - \log_{N\delta} Z(\XI)\Big|
	\le o_N(1) + o_A(1)
	\]
in the limit $N\to\infty$ followed by $A\to\infty$.

\begin{proof} We now interpolate between $Z(\XI)$ and $Z(u_A;\XI)$ by defining
	\begin{align*}
	Z_j
	&\equiv 
	\sum_{\sigma\in\set{-1,+1}^N}
	\exp\bigg\{\sum_{i=1}^j u_A\bigg( \f{(\XI^i,\sigma)}{N^{1/2}}\bigg)\bigg\}
	\cdot \prod_{i=j+1}^M
	U\bigg( \f{(g^i,\sigma)}{N^{1/2}}\bigg)\,, \\
	Z_{j,\circ}
	&\equiv 
	\sum_{\sigma\in\set{-1,+1}^N}
	\exp\bigg\{\sum_{i=1}^{j-1} u_A\bigg( \f{(\XI^i,\sigma)}{N^{1/2}}\bigg)\bigg\}
	\cdot \prod_{i=j+1}^M
	U\bigg( \f{(g^i,\sigma)}{N^{1/2}}\bigg)\,.
	\end{align*}
Then $Z_0= Z(\XI) \le Z_1 \le \ldots \le Z_M = Z(u_A;\XI)$, so 
	\[
	\f1N\Big| \log_{N\delta} Z(u_A;\XI) - \log_{N\delta} Z(\XI)\Big|
	=
	\f1N\sum_{j\le M}
	\Big( \log_{N\delta} Z_j- \log_{N\delta} Z_{j-1}\Big)
	\equiv \f1N\sum_{j\le M} Y_j\,.
	\]
Let $\bm{E}_{j,\circ}$ denote the event that
$Z_{j,\circ}\ge\exp(N\delta)$.
Note that $Z_{j,\circ}\ge \max\set{Z_j,Z_{j-1}}$, so $Y_j=0$ on the complement of $\bm{E}_{j,\circ}$. It follows that
$Y_j = \dot{y}_j-\ddot{y}_j$ where
	\begin{align*}
	\dot{y}_j
	&\equiv \Ind{\bm{E}_{j,\circ}}
	\Big( \log_{N\delta} Z_j- \log_{N\delta} Z_{j,\circ}\Big)
	=\Ind{\bm{E}_{j,\circ}}
	\log\max\bigg\{ \f{Z_j}{Z_{j,\circ}},
		\f{\exp(N\delta)}{Z_{j,\circ}} \bigg\}\,,\\
	\ddot{y}_j
	&\equiv\Ind{\bm{E}_{j,\circ}}
	 \Big( \log_{N\delta} Z_{j-1}- \log_{N\delta} Z_{j,\circ}\Big)
	= \Ind{\bm{E}_{j,\circ}}
	\log\max\bigg\{ \f{Z_{j-1}}{Z_{j,\circ}},
		\f{\exp(N\delta)}{Z_{j,\circ}} \bigg\}
	\,.
	\end{align*}
Writing $\<\cdot\>_{j,\circ}$ for expectation with respect to the Gibbs measure corresponding to $Z_{j,\circ}$, we have
	\begin{align*}
	\ddot{x}_j
	&\equiv \f{Z_{j-1}}{Z_{j,\circ}}
	=  \bigg\langle
	U\bigg( \f{(\XI^j,\sigma)}{N^{1/2}} \bigg)
	\bigg\rangle_{j,\circ}\,,\\
	\dot{x}_j
	&\equiv \f{Z_j}{Z_{j,\circ}}
	= \bigg\langle
	\exp u_A\bigg( \f{(\XI^j,\sigma)}{N^{1/2}} \bigg)
	\bigg\rangle_{j,\circ}
	= e^{-A}+(1-e^{-A}) \ddot{x}_j\,,
	\end{align*}
and clearly $0\le \ddot{x}_j \le \dot{x}_j \le 1$. Moreover we can bound
	\beq\label{e:U.A.ratio}
	1\le\f{\dot{x}_j}{\ddot{x}_j}
	=\f{e^{-A}+(1-e^{-A}) \ddot{x}_j}
		{\ddot{x}_j}
	=1+ \f{e^{-A}(1-\ddot{x}_j)}{\ddot{x}_j}
	\le 1+\f{1}{e^A\ddot{x}_j}\,.\eeq
Combining with \cite[Lem.~8.3.7]{MR3024566} (i.e., the bound \eqref{e:log.delta.bound} above) gives
	\begin{align}\nonumber
	0 \ge -Y_j &=\ddot{y}_j-\dot{y}_j
	\stackrel{\eqref{e:log.delta.bound}}{\ge}  \Ind{\bm{E}_{j,\circ}}
	\log \max\bigg\{ \f{\ddot{x}_j}{\dot{x}_j}, 
		\f{\exp(N\delta)}{Z_{j,\circ}}\bigg\} \\
	&\stackrel{\eqref{e:U.A.ratio}}{\ge}  \Ind{\bm{E}_{j,\circ}}
	\max\bigg\{
	-\log \bigg(1+\f1{e^A\ddot{x}_j}\bigg), 
	-N\log2
	\bigg\}
	\,.\label{e:log.lbd}
	\end{align}
Note $\log(1+x)\le x$ for all $x\ge0$,
and $\log(1+x)\le \log(2x)$ for $x\ge1$. 
Writing $\wmax\equiv N^{1/3}/C_\delta$, we hereafter assume $A\ge C_\delta$, and decompose
	\begin{align*}
	\textup{(I)}
	&\equiv
	\E\bigg[Y_j ; \ddot{x}_j \le
	\f1{\exp(\wmax)}\bigg]
	\stackrel{\eqref{e:log.lbd}}{\le} N\log2 \cdot
	\P\bigg( 
	\ddot{x}_j \le
	\f1{\exp(\wmax)};
	\bm{E}_{j,\circ}
	\bigg)\,,\\
	\textup{(II)}
	&\equiv 
	\E\bigg[Y_j ; 
	\f1{\exp(\wmax)}
	\le \ddot{x}_j \le \f1{e^A}
	\bigg]
	\stackrel{\eqref{e:log.lbd}}{\le}
	\E\bigg[
	\Ind{ \bm{E}_{j,\circ}}
	\log\bigg( \f{2}{e^A\ddot{x}_j}\bigg)
	 ;  \f1{\exp(\wmax)}
	\le \ddot{x}_j \le \f1{e^A}
	\bigg]\,,
	\\
	\textup{(III)}
	&\equiv\E\bigg[Y_j; 
	\f1{e^A} \le \ddot{x}_j \le 
	\f1{\exp(C_\delta)} 
	\bigg]
	\stackrel{\eqref{e:log.lbd}}{\le}
	\E\bigg[ \f{\Ind{ \bm{E}_{j,\circ}}}{e^A\ddot{x}_j};
	\f1{e^A} \le \ddot{x}_j \le 
	\f1{\exp(C_\delta)} 
	\bigg]\,.
	\end{align*}
If $U\equiv0$ there is nothing to prove, so we may assume $U(x)\ge\Ind{x\in[a,b]}$ for some $-\infty<a<b<\infty$. We can then apply Theorem~\ref{t:interval.general} to bound the above quantities:
	\begin{align*}
	\textup{(I)}
	&\le N\log2 \cdot
	\exp\bigg( -\f{\wmax}{C_\delta}\bigg)\,,\\
	\textup{(III)}
	&\le
	\f1{e^A}
	\int_{\exp(C_\delta)}^{\exp(A)} \f{du}{u^{1/C_\delta}} 
	\le
	\f1{e^A}
	\f{(e^A)^{1-1/C_\delta}}{1-1/C_\delta}
	\le \f{2}{e^{A/C_\delta}}\,.\end{align*}
Lastly, making a change of variables gives
	\[
	\P\bigg(\log\bigg( \f{2}{e^A \ddot{x}_j}\bigg) \ge t
		;\bm{E}_{j,\circ}\bigg)
	= \P\bigg(\f1{ \ddot{x}_j} \ge \f{e^{A+t}}{2}
		;\bm{E}_{j,\circ}\bigg)
	\le\bigg( \f{2}{e^{A+t}}\bigg)^{1/C_\delta}\,,
	\]
and integrating this tail bound over $t\ge0$ gives
	\[
	\textup{(II)}
	\le \int_0^\infty 
	\bigg( \f{2}{e^{A+t}}\bigg)^{1/C_\delta}\,dt
	\le \f{2 C_\delta}{e^{A/C_\delta}}\,.
	\]
Combining the bounds for $\textup{(I)}$, $\textup{(II)}$, and $\textup{(III)}$ gives $0\le \E Y_j \le o_A(1) + o_N(1)$, and the claim follows.
\end{proof}
\end{ppn}

\begin{proof}[\hypertarget{proof:t.univ.avg}{Proof of Theorem~\ref{t:univ.avg}}]
Follows by combining Proposition~\ref{p:univ.positive.temp}
with Proposition~\ref{p:pos.to.zero.temp}.
\end{proof}

The concentration of the free energy is addressed in the next section; see the \hyperlink{proof:t.conc}{proof of Theorem~\ref{t:conc}}.

\section{Concentration, sharp threshold sequence, and universality}
\label{s:sharp}

In this section we prove Theorems~\ref{t:conc}--\ref{t:univ}. Recall the abstract model \eqref{e:Z.Theta}. Since we mainly consider how the system behaves as $M$ varies, we will mostly drop $N$ from the notation, e.g.\ we will abbreviate $Z_M\equiv Z_{M,N}$. Let $\filt_a$ be the $\sigma$-field generated by the random functions $\Theta_1,\ldots,\Theta_a$. The following is a stronger version of Assumption~\ref{a:addone.weak}, which accommodates the perceptron models considered in this paper:

\begin{ass}\label{a:addone}
For the model \eqref{e:Z.Theta}, let $Z_{M+1}$ be the partition function that results from adding one more factor to $Z_M$. Suppose for all $\delta>0$ small enough that on the event $Z_M \ge \exp(N\delta)$ we have
	\[
	\P\bigg(
	\f{Z_{M+1}}{Z_M} \le \f1{\exp(w)} \,\bigg|\,
	\filt_M
	\bigg) \le \exp\bigg(-\f{w}{C_\delta}\bigg)
	\]
for all $C_\delta \le w \le \wmax$, where $C_\delta$ is a finite constant that depends only on the model and on $\delta$, and $\wmax$ can depend on the model as well as on $\delta$ and $N$.\end{ass}

In the case of the perceptron model, Assumption~\ref{a:addone} holds with the following parameters:
\begin{enumerate}[(a)]
\item For the half-space perceptron \eqref{e:Z.hspace} with gaussian disorder, Theorem~\ref{t:hspace.gaus} gives $\wmax=N/C_\delta$.

\item For the half-space perceptron \eqref{e:Z.hspace} with general disorder,
Theorem~\ref{t:hspace.general} gives $\wmax= N^{1/2}/C_\delta$.
\item For the $U$-perceptron \eqref{e:Z.U} 
with gaussian disorder,
Theorem~\ref{t:interval.gaus} gives $\wmax= N^{1/2}/C_\delta$.
\item For the $U$-perceptron \eqref{e:Z.U} with general disorder,
Theorem~\ref{t:interval.general} gives $\wmax= N^{1/3}/C_\delta$.
\end{enumerate}
(For both (c) and (d) above, we assume $U$ is $\set{0,1}$-valued with $U(x)\ge\Ind{x\in[a,b]}$ for some $-\infty<a<b<\infty$.) We also will show that the weaker Assumption~\ref{a:addone.weak} suffices for our main claims. This section is organized as follows:
\begin{itemize}
\item In \S\ref{ss:slow.decrease} we show that 
Assumption~\ref{a:addone.weak} or \ref{a:addone} implies
that the partition function of the model \eqref{e:Z.Theta} is unlikely 
to decrease very sharply after the addition of a small linear number of constraints.

\item In \S\ref{ss:conc.trunc} we show that
Assumption~\ref{a:addone.weak} or \ref{a:addone} implies
concentration of the truncated free energy of the model \eqref{e:Z.Theta}, leading to the 
\hyperlink{proof:t.conc}{proof of Theorem~\ref{t:conc}}.

\item In \S\ref{ss:sharp.univ} we combine the results obtained thus far to give the 
\hyperlink{p:t.sharpthreshold}{proof of Theorem~\ref{t:sharpthreshold}} (on the sharp threshold sequence) and finally the
\hyperlink{proof:t.univ}{proof of Theorem~\ref{t:univ}} (on universality).

\end{itemize}

\subsection{Slow decrease of partition function with new constraints}
\label{ss:slow.decrease}

In this subsection we prove 
Proposition~\ref{p:slowdecrease}, which says that under 
Assumption~\ref{a:addone}, the partition function of the model \eqref{e:Z.Theta} is unlikely 
to decrease very sharply after adding a small linear number of constraints. We also prove
Proposition~\ref{p:slowdecrease.weak}
which gives a similar but weaker estimate under
Assumption~\ref{a:addone.weak}.

\begin{ppn}\label{p:slowdecrease}
For the model \eqref{e:Z.Theta}, under Assumption~\ref{a:addone}, 
for all $\delta>0$ there exists $\rho_\delta>0$ such that for all $0\le\rho\le\rho_\delta$ we have
	\[
	\P\bigg( Z_{M+N\rho} < \exp(N\delta)
		\,\bigg|\, Z_M \ge \exp(2N\delta)\bigg)
	\le N\rho \exp\bigg(-\f{\wmax}{C_\delta}\bigg)
	+ \exp\bigg(-\f{N\delta}{4C_\delta}\bigg)\,.\]
In particular, if $\wmax$ grows faster than $\log N$, then this is $o_N(1)$.

\begin{proof}For $0\le\ell\le N\rho$ let us abbreviate $X_\ell=Z_{M+\ell}$ and $\GG_\ell\equiv\filt_{M+\ell}$. It follows by Assumption~\ref{a:addone} that on the event $X_\ell\ge\exp(N\delta)$, we have the bound
	\beq\label{e:assump.rewrite}
	\P\bigg(
	\f{X_{\ell+1}}{X_\ell} \le \f1{\exp(w)}\,\bigg|\,\GG_\ell\bigg)
	\le \exp\bigg(-\f{w}{C_\delta}\bigg)\,,
	\eeq
provided $C_\delta \le w\le \wmax$.
To apply this bound, we define a stopping time $0\le\ZETA\le N\rho$ by
	\beq\label{e:slowdec.ZETA}
	\ZETA\equiv\min\bigg\{\ell\ge0: \ell=N\rho
	\textup{ or }
	X_\ell <\exp(N\delta)\bigg\}\,.
	\eeq
Applying \eqref{e:assump.rewrite} with $w=\wmax$ gives, on the event $X_0\ge\exp(2N\delta)$,
	\beq\label{e:slowdec.unionbound}
	\textup{(I)}
	\equiv \P\bigg(
	\min_{\ell < \ZETA}
	\f{X_{\ell+1}}{X_\ell} \le \f1{\exp(\wmax)}
	\,\bigg|\,\GG_0\bigg)
	\le N\rho \exp\bigg( - \f{\wmax}{C_\delta}\bigg)\,.
	\eeq
Next define the truncated random variables
	\[V_{\ell+1} \equiv 
	\min\bigg\{ \f{X_\ell}{X_{\ell+1}},\exp(\wmax)\bigg\}
	=
	\min\bigg\{ \f{Z_{M+\ell}}{Z_{M+\ell+1}},\exp(\wmax)\bigg\}
	\,.\]
On the event $X_\ell \ge\exp(N\delta)$, \eqref{e:assump.rewrite} implies the tail bound
	\[
	\P\bigg( (V_{\ell+1})^{1/(2C_\delta)} \ge y
	\,\bigg|\, \GG_\ell\bigg)
	\le \f1{y^2}\,,
	\]
for all $y\ge \exp(C_\delta/2)$. Integrating over $y$ gives the moment bound
	\beq\label{e:moment.bd}
	\E\bigg[ (V_{\ell+1})^{1/(2C_\delta)} 
	\,\bigg|\, \GG_\ell\bigg]
	\le \exp\bigg( \f{C_\delta}{2}\bigg)
	 + \int_{\exp(C_\delta/2)}^\infty \f{dy}{y^2}
	\le \exp(C_\delta)\,.
	\eeq
Applying Markov's inequality gives, on the event $Z_M=X_0\ge\exp(2N\delta)$,
	\[
	\textup{(II)}
	\equiv
	\P\bigg(
	\prod_{\ell\le\ZETA} V_\ell \ge \exp(N\delta)
	\,\bigg|\,\GG_0
	\bigg)
	\le 
	\f{\exp(N\rho C_\delta)}{\exp(N\delta/(2C_\delta))}
	\le \exp\bigg(-\f{N\delta}{4C_\delta}\bigg)\,,
	\]
where the last bound holds by taking $\rho\le \rho_\delta$ small enough. Thus, on the event $Z_M=X_0\ge\exp(2N\delta)$, we have
	\[
	\P\bigg( X_\ZETA < \exp(N\delta)
		\,\bigg|\, \GG_0 \bigg)
	\le \textup{(I)}+\textup{(II)}
	\le
	N\rho \exp\bigg(-\f{\wmax}{C_\delta}\bigg)
	+ \exp\bigg(-\f{N\delta}{4C_\delta}\bigg)\,.
	\]
On the complementary event $X_\ZETA\ge \exp(N\delta)$,
it follows from the definition \eqref{e:slowdec.ZETA} that we must have $\ZETA=N\rho$, which concludes the proof.
\end{proof}
\end{ppn}

\begin{ppn}\label{p:slowdecrease.weak}
For the model \eqref{e:Z.Theta}, under Assumption~\ref{a:addone.weak}, 
for all $\delta>0$ there exists $\rho_\delta>0$ such that for all $0\le\rho\le\rho_\delta$ we have 
	\[
	\P\bigg( Z_{M+N\rho} < \exp(N\delta)
		\,\bigg|\, Z_M \ge \exp(2N\delta)\bigg)
	\le
	N\rho f_\delta(\wmax)+\f{4\rho C_{\delta,2}}{N\delta^2}\,.\]
In particular, if $f_\delta(\wmax)$ tends to zero superpolynomially in $N$, then this is $o_N(1)$.

\begin{proof}
Following the first few steps of the proof of Proposition~\ref{p:slowdecrease},
with $X_\ell\equiv Z_{M+\ell}$,
$\GG_\ell\equiv\filt_{M+\ell}$, and the same stopping time $\ZETA$ as defined by \eqref{e:slowdec.ZETA},
yields that on the event $X_0\ge\exp(2N\delta)$ we have
	\[
	\textup{(I)}
	\equiv \P\bigg(
	\min_{\ell < \ZETA}
	\f{X_{\ell+1}}{X_\ell} \le \f1{\exp(\wmax)}
	\,\bigg|\,\GG_0
	\bigg)
	\le N\rho f_\delta(\wmax)\,.
	\]
(instead of \eqref{e:slowdec.unionbound}). Next define the truncated random variables
	\[V_{\ell+1} \equiv \min\bigg\{ 
	\f{X_\ell}{X_{\ell+1}}, \exp(\wmax)\bigg\}\,,\]
and note that $V_{\ell+1}\ge1$. Assumption~\ref{a:addone.weak} implies
that on the event $X_\ell\ge \exp(N\delta)$ we have
	\begin{align*}
	\E\bigg(\log V_{\ell+1}
	\,\bigg|\,\GG_\ell\bigg)
	&\le C_\delta+ \int_{C_\delta}^\infty f_\delta(w)\,dw
	\le C_{\delta,1}\,. \\
	\E\bigg((\log V_{\ell+1})^2
	\,\bigg|\,\GG_\ell\bigg)
	&\le (C_\delta)^2
	+ \int_{C_\delta}^\infty 2w f_\delta(w)\,dw
	\le C_{\delta,2}\,.
	\end{align*}
Now consider the Doob martingale decomposition 
	\[
	\sum_{k=1}^\ell \log V_k
	= 
	\sum_{k=1}^\ell
	\Big(\log V_k-\E(\log V_k\,|\,\GG_{k-1})
	\Big)
	+\sum_{k=1}^\ell\E(\log V_k\,|\,\GG_{k-1})
	\equiv 
	M_\ell + A_\ell
	\]
It follows from the preceding bounds that  $0\le A_\ZETA \le N\rho C_{\delta,1}$ and $\E[(M_\ZETA)^2] \le N\rho C_{\delta,2}$.
We can take $\rho$ small enough to guarantee that $\rho C_{\delta,1} <\delta/2$. It follows using Chebychev's inequality that, on the event $X_0\ge\exp(2N\delta)$,
	\[
	\textup{(II)}
	\equiv \P\bigg(\sum_{k\le\ZETA} \log V_k
	\ge N\delta \,\bigg|\, \GG_0\bigg)
	\le\P\bigg(M_\ZETA \ge N\delta - N\rho C_{\delta,1}
	\ge \f{N\delta}{2} \,\bigg|\, \GG_0\bigg)
	\le \f{4\rho C_{\delta,2}}{N\delta^2}\,.
	\]
Combining the above bounds for $\textup{(I)}$ and $\textup{(II)}$ gives,
on the event $X_0\ge\exp(2N\delta)$,
	\[
	\P\bigg( X_\ZETA < \exp(N\delta)
		\,\bigg|\, \GG_0\bigg)
	\le \textup{(I)}+\textup{(II)}
	\le
	N\rho f_\delta(\wmax)+\f{4\rho C_{\delta,2}}{N\delta^2}\,.
	\]
On the complementary event $X_\ZETA\ge \exp(N\delta)$,
it follows from the definition \eqref{e:slowdec.ZETA} that we must have $\ZETA=N\rho$, and this concludes the proof.
\end{proof}
\end{ppn}

\subsection{Concentration of truncated free energy}
\label{ss:conc.trunc}

 The first result below, 
Proposition~\ref{p:truncated.conc}, is fairly similar to the result of 
\cite[Propn.~9.2.6]{MR3024566}, and makes use of Assumption~\ref{a:addone}. The next result, Proposition~\ref{p:truncated.conc.weak}, is a similar but weaker estimate that uses only Assumption~\ref{a:addone.weak}. We conclude the subsection with the \hypertarget{proof:t.conc}{proof of Theorem~\ref{t:conc}}.

\begin{ppn}\label{p:truncated.conc}
For the model \eqref{e:Z.Theta}, under Assumption~\ref{a:addone}, 
for all $\delta>0$ we have
	\[
	\P\bigg( \Big|\log_{N\delta} Z
	-\E\log_{N\delta} Z\Big|\ge 2Nt\bigg)
	\le
	\f{M\log 2}{t} \exp\bigg(-\f{\wmax}{C_\delta}\bigg)
	+\exp\bigg( - \min\bigg\{ \f{Nt^2}{4 (A_\delta)^2},
	\f{Nt}{2A_\delta} \bigg\}\bigg)\,,
	\]
where $A_\delta$ is a large constant which depends only on $C_\delta$  (where $C_\delta$ is the constant appearing in Assumption~\ref{a:addone}). In particular, if $\wmax$ grows faster than $\log N$, then the above is 
$o_N(1)$ provided $t\gg 1/N^{1/2}$.

\begin{proof}
Recall that $\filt_k$ denotes the $\sigma$-field generated by the random functions $\Theta_1,\ldots,\Theta_k$. Consider the martingale decomposition
	\[
	\log_{N\delta}Z-\E(\log_{N\delta}Z)
	= \sum_{k\le M}
	\bigg\{
	\E(\log_{N\delta}Z\,|\,\filt_k)
	-\E(\log_{N\delta}Z\,|\,\filt_{k-1})\bigg\}
	\equiv \sum_{k\le M} Y_k\,,
	\]
Let $Z_{k,\circ}$ be the partition function without the $k$-th factor:
	\[
	Z_{k,\circ}\equiv \sum_{\sigma\in\set{-1,+1}^N}
	\prod_{\ell\le M, \ell \ne k} \Theta_\ell(\sigma)\,.\]
Since we assumed the $\Theta$ functions are $\set{0,1}$-valued, we have $Z \le Z_{k,\circ}$, so if $Z_{k,\circ} \le \exp(N\delta)$ then we must have $\log_{N\delta}Z=\log_{N\delta}Z_{k,\circ}$. It follows that
	\[
	L_k
	\equiv \log_{N\delta}Z-\log_{N\delta}Z_{k,\circ}
	=\bigg( \log_{N\delta}Z-\log Z_{k,\circ}\bigg)
	\Ind{Z_{k,\circ}\ge \exp(N\delta)}\,.
	\]
Let $\E_k$ denote expectation over $\Theta_k$ only,
and let $\E^k$ denote expectation over all the $(\Theta_\ell)_{\ell\ge k}$. We can then rewrite
	\[
	Y_k
	=
	\E\bigg(\log_{N\delta}Z-\log_{N\delta}Z_{k,\circ}\,\bigg|\,\filt_k\bigg)
	-\E\bigg(\log_{N\delta}Z-\log_{N\delta}Z_{k,\circ}\,\bigg|\,\filt_{k-1}\bigg)
	= \E^{k-1}(L_k - \E_k L_k)\,.
	\]
For comparison, let us also define
$y_k\equiv\E^{k-1}(\ell_k - \E_k \ell_k)$ where
	\[
	\ell_k
	\equiv
	\bigg( \log_{N\delta}Z-\log Z_{k,\circ}\bigg)
	\mathbf{1}\bigg\{Z_{k,\circ}\ge \exp(N\delta) ; 
	\f{Z}{Z_{k,\circ}} \ge \f1{\exp(\wmax)}\bigg\}
	\]
Since $-N\log2\le \log_{N\delta}Z-\log_{N\delta} Z_{k,\circ} \le0$, we can  use Assumption~\ref{a:addone} to bound
	\begin{align}\nonumber
	\textup{(I)}
	&\equiv
	\P\bigg(\sum_{k\le M} |Y_k-y_k| \ge Nt\bigg)
	\le \f1{Nt} \E\bigg(\sum_{k\le M} |Y_k-y_k|\bigg) \\
	& \le \f{MN\log2}{Nt}
		\P\bigg( Z_{k,\circ}\ge \exp(N\delta);
		\f{Z}{Z_{k,\circ}} \ge \f1{\exp(\wmax)} \bigg)
	\le \f{MN\log 2}{Nt}
		\exp\bigg(-\f{\wmax}{C_\delta}\bigg)\,.
	\label{e:largedevs}
	\end{align}
We will bound, for small enough $\theta$, the exponential moment
	\begin{align*}
	\E\bigg(\exp(\theta|y_k|)\,\bigg|\,\filt_{k-1}\bigg)
	&= \E^{k-1}\E_k \exp\bigg\{ \theta 
		\Big|\E^{k-1}(\ell_k - \E_k \ell_k)\Big| \bigg\}
	\le
	\E^{k-1}\E_k \exp \Big(\theta |\ell_k- \E_k \ell_k|\Big) \\
	&\le
	\E^{k-1}\bigg[
	\Big(\E_k \exp (\theta|\ell_k|)\Big)
	\cdot \exp(\theta|\E_k \ell_k|)
	\bigg]
	\le \E^{k-1}\bigg[\Big(\E_k \exp (\theta|\ell_k|)\Big)^2\bigg]\,.
	\end{align*}
We now proceed to bound $\E_k \exp (\theta|\ell_k|)$: for $\theta=1/(2C_\delta)\ge0$, Assumption~\ref{a:addone} gives
	\[\E_k \exp(\theta|\ell_k|)
	\le
	\E_k\bigg[
	\bigg(\f{Z_{k,\circ}}{Z}\bigg)^{1/(2C_\delta)}
	; Z_{k,\circ} \ge \exp(N\delta),
	\f{Z}{Z_{k,\circ}} \le \f1{\exp(\wmax)}\bigg]
	\le \exp(C_\delta)\,,
	\]
by the same calculation as in \eqref{e:moment.bd}. It follows using Jensen's inequality that for $A_\delta$ a large enough constant (depending only on $C_\delta$) we have 
	\[\E_k \exp\bigg(\f{|\ell_k|}{A_\delta}\bigg)
	\le 
	\bigg( \E_k \exp\bigg(\f{|\ell_k|}{2C_\delta}\bigg)
	\bigg)^{2C_\delta/A_\delta}
	\le \exp\bigg(\f{2(C_\delta)^2}{A_\delta}\bigg)
	\le2\,,\]
where the last inequality holds by choosing $A_\delta$ large enough (depending only on $C_\delta$). It follows by the martingale Bernstein inequality (see \cite[Thm.~A.6.1]{MR3024566}) that
	\beq\label{e:bernstein}
	\textup{(II)}
	\equiv\P\bigg(\bigg|\sum_{k\le M} y_k\bigg|\ge Nt\bigg)
	\le
	\exp\bigg( - \min\bigg\{ \f{Nt^2}{4 (A_\delta)^2},
	\f{Nt}{2A_\delta} \bigg\}\bigg)\,.\eeq
The claim follows by combining \eqref{e:largedevs} and \eqref{e:bernstein}.
\end{proof}
\end{ppn}

\begin{ppn}\label{p:truncated.conc.weak}
For the model \eqref{e:Z.Theta}, under Assumption~\ref{a:addone.weak}, 
for all $\delta>0$ we have
    \[
	\P\bigg( \Big|\log_{N\delta} Z
	-\E\log_{N\delta} Z\Big|\ge 2Nt\bigg)
	\le
	\f{M\log 2}{t} f_\delta(\wmax)
	+\f{MC_{\delta,2}}{(Nt)^2}\,.
	\]
In particular, if $f_\delta(\wmax)$ tends to zero superpolynomially in $N$, then this is $o_N(1)$ provided $t\gg 1/N^{1/2}$.

\begin{proof}
Following the first few steps of the proof of Proposition~\ref{p:truncated.conc} gives (instead of \eqref{e:largedevs})
	\[
	\textup{(I)}
	\equiv
	\P\bigg(\sum_{a\le M} |Y_a-y_a| \ge Nt\bigg)
	\le \f1{Nt} \E\bigg(\sum_{a\le M} |Y_a-y_a|\bigg)
	\le \f{MN\log 2 }{Nt} f_\delta(\wmax)\,.\]
Recall that $\E_k$ denotes expectation over $\Theta_a$ only,
while $\E^k$ denotes expectation over all $(\Theta_\ell)_{\ell\ge k}$. Then
	\[
	\E_a[ (y_a)^2]
	\le \E_a[ (\ell_a)^2]
	\le \E_a\bigg[
	\bigg(\log \f{Z}{Z_{a,\circ}}\bigg)^2
	; \f{Z}{Z_{a,\circ}} \ge \f1{\exp(\wmax)}\bigg]
	\le C_{\delta,2}\,,
	\]
by Assumption~\ref{a:addone.weak}. It follows using Chebychev's inequality that
	\[\textup{(II)}
	\equiv
	\P\bigg( \bigg|\sum_{a\le M} y_a\bigg| \ge Nt\bigg)
	\le \f{MC_{\delta,2}}{(Nt)^2}\,.
	\]
The claim follows by combining the bounds on $\textup{(I)}$ and $\textup{(II)}$.
\end{proof}
\end{ppn}

\begin{proof}[\hypertarget{proof:t.conc}{Proof of Theorem~\ref{t:conc}}]
Follows from Proposition~\ref{p:truncated.conc.weak}. In the perceptron model, a sharper concentration result can be obtained using Proposition~\ref{p:truncated.conc}, where Assumption~\ref{a:addone} is satisfied by Theorems~\ref{t:hspace.gaus}--\ref{t:interval.general}.
\end{proof}

\subsection{Sharp threshold sequence and universality} 
\label{ss:sharp.univ} In this subsection we conclude with the
proofs of Theorems~\ref{t:sharpthreshold} and \ref{t:univ}. We summarize the results of the preceding subsections with the following assumptions:

\begin{ass}\label{a:slow.decay}
For all $\delta>0$ there exists $\rho_\delta>0$ such that for all $0\le\rho\le\rho_\delta$ we have, for all $M\ge0$,
	\[
	\P\bigg(Z_M \ge \exp(2N\delta), Z_{M+N\rho} < \exp(N\delta)\bigg)
	\le o_N(1)\,.\]
(Under Assumptions~\ref{a:addone} or \ref{a:addone.weak}, this estimate is implied by Propositions~\ref{p:slowdecrease} and \ref{p:slowdecrease.weak}.)
\end{ass}

\begin{ass}\label{a:trunc.conc}
For all $\delta>0$ and any constant $t>0$,
	\[
	\P\bigg(\Big|\log_{N\delta}Z-\E\log_{N\delta}Z\Big|\ge Nt\bigg)
	\le o_N(1)\,.
	\]
(Under Assumptions~\ref{a:addone} or \ref{a:addone.weak}, this estimate is implied by Propositions~\ref{p:truncated.conc} and \ref{p:truncated.conc.weak}.)
\end{ass}

\begin{ppn}\label{p:sharp.threshold}
For the model \eqref{e:Z.Theta}, suppose $\sup\{\E\Theta_a(x):x\in\set{-1,+1}^N\}\le \exp(-c)$ for a positive constant $c$. Assumption~\ref{a:trunc.conc} implies a sharp threshold sequence: that is to say, there is a sequence $\alpha_N$ such that $\P(Z_{N\alpha}>0)$ transitions from $1-o_N(1)$ to $o_N(1)$ in an $o_N(1)$ window around $\alpha_N$. Assumption~\ref{a:slow.decay} further implies $\alpha_N\asymp1$ in the limit $N\to\infty$.

\begin{proof} 
Note that $Z_0=2^N$. The bound on $\E\Theta_a$ implies $\E Z_M \le 2^N\exp(-Mc)$, so $\P(Z_{N\alpha}>0)=o_N(1)$ as soon as $\alpha > (\log2)/c$.
Thus $\P(Z_{N\alpha}>0)$ transitions from $1-o_N(1)$ to $o_N(1)$ as $\alpha$ increases from $0$ to $2(\log 2)/c$. Now suppose for the sake of contradiction that there is not a sharp threshold sequence: this means that there exists arbitrarily large $N$ such that we have
	\beq\label{e:no.sharp.threshold}
	\epsilon
	\le
	\inf\bigg\{
	\P(Z_{N\alpha}>0) : \alpha\in[\alpha_1,\alpha_2]
	\bigg\}
	\le
	\sup\bigg\{
	\P(Z_{N\alpha}>0) : \alpha\in[\alpha_1,\alpha_2]
	\bigg\}\le1-\epsilon\,,\eeq
where the $\alpha_i$ can depend on $N$, but $\alpha_2-\alpha_1$ stays bounded away from zero as $N\to\infty$. Then
	\begin{align*}
	\P\bigg(Z_{N\alpha_1}
	\le \exp\bigg(\f{N(\alpha_2-\alpha_1)c}{2}\bigg)\bigg)
	&\le \P(Z_{N\alpha_2}=0)
	+ \P\bigg(Z_{N\alpha_1}
		\le \exp\bigg(\f{N(\alpha_2-\alpha_1)c}{2}\bigg),
			Z_{N\alpha_2}>0\bigg) \\
	&\le \P(Z_{N\alpha_2}=0) 
		+ \exp\bigg(-\f{N(\alpha_2-\alpha_1)c}{2}\bigg)
	\le 1-\epsilon + o_N(1)\,.
	\end{align*}
In the above, the second-to-last step is by the assumption on $\E\Theta$,
and the last step is by \eqref{e:no.sharp.threshold}. Rearranging the above inequality gives
	\[\epsilon-o_N(1)
	\le \P\bigg(Z_{N\alpha_1}
	> \exp\bigg(\f{N(\alpha_2-\alpha_1)c}{2}\bigg)\bigg)
	\le \P(Z_{N\alpha_1}>0) \le 1-\epsilon\,.
	\]
It follows that for any $\delta < (\alpha_2-\alpha_1)c/4$, the quantity $\log_{N\delta} Z_{N\alpha_1}$ is not well-concentrated:
with probability at least $\epsilon-o_N(1)$ it equals $N\delta$, but 
with probability at least $\epsilon-o_N(1)$ it exceeds $2N\delta$. This contradicts Assumption~\ref{a:trunc.conc}. This implies that 
$\P(Z_{N\alpha}>0)$ transitions from $1-o_N(1)$ to $o_N(1)$ in an $o_N(1)$ window around a sharp threshold sequence $\alpha_N$. Finally, Assumption~\ref{a:slow.decay} implies that $\alpha_N$ stays bounded away from zero in the limit $N\to\infty$.
\end{proof}
\end{ppn}

\begin{proof}[\hypertarget{p:t.sharpthreshold}{Proof of Theorem~\ref{t:sharpthreshold}}]
Follows from Proposition~\ref{p:sharp.threshold}, where 
Assumptions~\ref{a:slow.decay} and \ref{a:trunc.conc} 
are satisfied by Propositions \ref{p:slowdecrease.weak} and \ref{p:truncated.conc.weak}.
\end{proof}

\begin{proof}[\hypertarget{proof:t.univ}{Proof of Theorem~\ref{t:univ}}]
For the $U$-perceptron \eqref{e:Z.U}, let $Z(\XI)$ be the partition function with general (subgaussian) disorder. It follows from the proof of Theorem~\ref{t:sharpthreshold} that there is a sharp threshold sequence $\alpha_{N,\XI}$ (potentially depending on the disorder), and the free energy $\log Z_{N,N\alpha}(\XI)$ is exponentially large for $\alpha<\alpha_{N,\XI}$. Suppose the threshold sequence depends nontrivially on the disorder, meaning that on a subsequence $N\to\infty$ we have
$\alpha_{N,\XI}-\alpha_{N,g}\ge\epsilon$ or
$\alpha_{N,g}-\alpha_{N,\XI}\ge\epsilon$.
Without loss suppose $\alpha_{N,\XI}-\alpha_{N,g}\ge\epsilon$.
Then, for $\alpha_{N,g}<\alpha<\alpha_{N,\XI}$,
 $Z_{N,N\alpha}(g)$ will be exponentially large while 
 $Z_{N,N\alpha}(\XI)$ is zero with high probability.
By the concentration result Theorem~\ref{t:conc}, this will yield a contradiction to Theorem~\ref{t:univ.avg}. It follows $Z(g)$ and $Z(\XI)$ share the same sharp threshold sequence $\alpha_N$. Moreover, since the partition function must be exponentially large for $\alpha<\alpha_N$, it follows from Theorem~\ref{t:conc} and 
Theorem~\ref{t:univ.avg} that $N^{-1}\log Z(\XI)$ and $N^{-1}\log Z(g)$
are $o_N(1)$-close with high probability for all $\alpha<\alpha_N$.
\end{proof}

{\raggedright
\bibliographystyle{alphaabbr}
\bibliography{prefs}
}

\end{document}